\documentclass[11pt,reqno]{amsart}
\usepackage{color}
\usepackage{bm}
\usepackage{marginnote}
\usepackage{amsmath,amssymb}
\usepackage{verbatim}

\catcode`\@=11 \@addtoreset{equation}{section}

\catcode`\@=12

\newtheorem{Theorem}{Theorem}[section]
\newtheorem{Proposition}{Proposition}[section]
\newtheorem{Lemma}{Lemma}[section]
\newtheorem{Corollary}{Corollary}[section]
\newtheorem{Definition}{Definition}[section]
\newtheorem{Remark}{Remark}[section]

\newcommand{\Div}	{{\rm div}_x}
\newcommand{\R} {\mathbb{R}}
\newcommand{\N} {\mathbb{N}}
\newcommand{\D} {\mathcal{D}}
\newcommand{\B} {\mathcal{B}}
\newcommand{\A} {\mathcal{A}}
\newcommand{\Nabla} {\nabla_x}
\newcommand{\vu}{\mathbf{u}}
\newcommand{\vm}{\mathbf{m}}
\newcommand{\vn}{\mathbf{n}}
\newcommand{\dx}{{\rm\, d}x}
\newcommand{\dy}{{\rm\, d}y}
\newcommand{\dt} {{\rm\, d}t}
\newcommand{\dxdt} {{\rm\,d}x{\rm d}t}
\newcommand{\vvarphi} {\mathbf\varphi}
\newcommand{\ep}{\varepsilon}
\newcommand{\weak}			{\rightharpoonup}

\newcommand{\vc}{\mathbf}
\newcommand{\Grad}{\nabla}
\newcommand{\intO}[1]{\int_{\Omega} #1 \ \dx}
\newcommand{\intTO}[1]{\int_0^T\!\!\!\! \int_{\Omega} #1 \ \dxdt}

\newcommand{\lr}[1]{\left( #1 \right)}
\newcommand{\Ov}[1]{\overline{#1}}
\newcommand{\vr}{\varrho}
\newcommand{\eq}[1]{\begin{equation}
\begin{split}
#1
\end{split}
\end{equation}}
\newcommand{\eqh}[1]{\begin{equation*}
\begin{split}
#1
\end{split}
\end{equation*}}

\definecolor{grey}{rgb}{0.85,0.85,0.85}
\date{}

\begin{document}

%%%%%%%%%%%%%%%%%%%%%%%%%%%%%%%%

\title[Long-time asymptotics for hydrodynamics with balance of forces]{On long-time asymptotics for viscous hydrodynamic models of collective behavior with damping and nonlocal interactions}

\author{Jos\'{e} A. Carrillo}
\address[Jos\'{e} A. Carrillo]{\newline Department of Mathematics, Imperial College London, London SW7 2AZ, United Kingdom}
\email{carrillo@imperial.ac.uk}

\author{Aneta Wr\'oblewska-Kami\'nska}
\address[Aneta Wr\'oblewska-Kami\'nska]{\newline Department of Mathematics, Imperial College London, London SW7 2AZ, United Kingdom and Institute of Mathematics, Polish Academy of Sciences, \'Sniadeckich 8, 00-656 Warszawa, Poland}
\email{a.wroblewska-kaminska@imperial.ac.uk,awrob@impan.pl}

\author{Ewelina Zatorska}
\address[Ewelina Zatorska]{\newline Department of Mathematics, University College London,  Gower Street, London WC1E 6BT, United Kingdom}
\email{e.zatorska@ucl.ac.uk}

\maketitle

\begin{abstract}
Hydrodynamic systems arising in swarming modelling include nonlocal forces in the form of attractive-repulsive potentials as well as pressure terms modelling strong local repulsion.  We focus on the case where there is a balance between nonlocal attraction and local pressure in presence of confinement in the whole space. Under suitable assumptions on the potentials and the pressure functions, we show the global existence of weak solutions for the hydrodynamic model with viscosity and linear damping. By introducing linear damping in the system, we ensure the existence and uniqueness of stationary solutions with compactly supported density, fixed mass and center of mass. The associated velocity field is zero in the support of the density.  Moreover, we show that global weak solutions converge for large times to the set of these stationary solutions in a suitable sense. In particular cases, we can identify the limiting density uniquely as the global minimizer of the free energy with the right mass and center of mass.
\end{abstract}

\medskip

\medskip

\section{Introduction}\label{i}
Continuum hydrodynamic descriptions for collective behavior of particles/agents are a very useful tool in mathematical biology to efficiently model the behavior of large populations of cells moving due to interactions produced by adhesion or chemical cues, and in large groups of animals via visual or sensory interactions. These models have been derived either phenomenologically, as in \cite{TT,GACCDGSPB}, or by the methods of kinetic theory, see \cite{CDMBC,CDP,CCR,dm2008,review,CKMT,dfl2013,dfl2015,AP,KT,CCP} and the references therein. Particles are assumed to interact nonlocally via attractive and repulsive forces modelling a range of these effects such as cell adhesion, chemotaxis interaction or volume constraints. Nonlinear pressure has been used as a model for volume size effects in cell or animal populations  \cite{TT,GACCDGSPB,HPVolume,CaCa06} as it can be seen as a very localized repulsive interaction.
In this work, we are interested in qualitative properties of the following hydrodynamic system  
	\begin{equation}\label{main_system}
	\begin{split}
	\partial_t \varrho  + \Div ( \varrho \vu) 
	& = 0 \, \\
	\partial_t \varrho \vu + \Div ( \varrho \vu \otimes \vu ) + a \Nabla  \varrho^m 
	&=  \mu \Delta \vu + \nabla_x (\lambda + \mu) \Div \vu- (\Nabla K \ast \varrho) \varrho -\varrho  \Nabla \Phi - \varrho \vu
	 \,
	\end{split}
	\end{equation}	
considered in $(0,T) \times \Omega$, where the unknowns are: $\varrho(t): \Omega \to \R$ for $t\ge 0$ denoting the density and $\vu(t) : \Omega \to \R^3$ for $t\ge 0$ denoting the velocity field. Moreover, the interaction potential $K: \R^3 \to \R$ encodes the nonlocal interactions (repulsive or attractive), $\Phi: \R^3 \to \R$ is a confinement potential that maybe present or not, and the physical constants involved in the viscosity term and the nonlinear pressure satisfies the following assumptions
	\begin{equation}\label{lambda_mu_m}
	\mu >0, \quad  \lambda+\frac{2}{3} \mu \geq 0, \quad a>0, \quad m>\frac{3}{2}\,.
	\end{equation}
The restriction for the nonlinear pressure exponent $m$ allows for compactness properties for solutions and approximate solutions. 

The system \eqref{main_system} will be either considered in the whole space $\Omega=\R^3$ with a confinement potential $\Phi$ and the right functional setting or in a bounded smooth domain supplemented with the Dirichlet boundary condition:
\begin{equation}\label{bc_intro}
	\vu |_{\partial \Omega} = 0,
\end{equation}
with or without the confinement potential. In that case the convolution in the nonlocal term $\Nabla K \ast \varrho$ is defined by extending the density to the whole space by zero outside its domain of definition. 

Our problem is supplemented with initial data $( \varrho_0, \vm_0 ) $ such that 
	\begin{equation}\label{initial_data}
	\varrho(0,x)=\varrho_0  \in L^1_+(\Omega)\cap L^m(\Omega), \quad \varrho \vu(0,x) = \vm_0 \in L^1(\Omega)\cap L^{\frac{6}{5}}(\Omega)  
	\end{equation}
and the following compatibility condition is supposed to be satisfied
	\begin{equation}\label{com_con}
	\vm_0= 0 \mbox{ whenever } \varrho_0=0, \frac{|\vm_0|^2}{\varrho_0} \in L^1(\Omega).
	\end{equation}

Stationary solutions to  \eqref{main_system} with $\vu=0$ in the support of $\vr$ are intimately related to aggregation-diffusion equations of the type
\begin{equation}\label{aggdif}
\partial_t \varrho  + \Div ( \varrho \bar\vu) = 0 \qquad \mbox{with } \bar\vu =
- a \Nabla \varrho^m - (\Nabla K \ast \varrho) \varrho -\varrho  \Nabla \Phi \,.
\end{equation}
Let us first point out that \eqref{aggdif} can be obtained through a relaxation limit from \eqref{main_system} as in \cite{LT} for zero viscosity coefficients. On the other hand, if there is a solution $\varrho_s$ of the problem
\begin{equation}\label{balanceforces}
a \Nabla \varrho^m_s + (\Nabla K \ast \varrho_s) \varrho_s + \varrho_s  \Nabla \Phi=0\,,
\end{equation}
in $\Omega$, then the pair $(\varrho_s,\vu_s=0)$ is a stationary solution of the hydrodynamic system \eqref{main_system}. The connection between these two macroscopic models \eqref{main_system} and \eqref{aggdif} is deeply rooted in their variational structure. Both the hydrodynamic system \eqref{main_system} and the aggregation-diffusion equation \eqref{aggdif} dissipate the total free energy defined as
\begin{equation}\label{energytot}
E(t) =  \int_{\Omega} \frac{1}{2}\varrho |\vu|^2 \dx + \mathcal{F}[\rho] \qquad \mbox{with} \qquad \mathcal{F}[\rho] =  \int_{\Omega} \left[ \frac{a}{m-1} \varrho^m + \frac{1}{2} (K \ast\varrho)\varrho  +  \varrho \Phi \right]\dx\,.
\end{equation}
Moreover, the dissipation of the free energy for the hydrodynamic system \eqref{main_system}, formally given by
\begin{equation}\label{endis}
\frac{d}{dt} E(t) = -\int_{\Omega}\left[ \mu |\nabla \vu |^2 + (\lambda+ \mu) | {\rm div}\vu|^2 +   \varrho |\vu|^2  \right] \dx \,,
\end{equation}
vanishes for zero velocity field due to the boundary conditions. Therefore, the stationary solutions for the hydrodynamic system \eqref{main_system} have a density satisfying \eqref{balanceforces}. 

In fact, finding conditions leading to a balance between repulsion, modelled by nonlinear diffusion, and aggregation, modelled by nonlocal interactions, has been very popular in the last 10 years due to its importance in mathematical biology and other applications for particular potentials in gravitational mechanics. Finding stationary densities that satisfy the balance of forces \eqref{balanceforces} is a very challenging question by itself. It is related to finding minimizers of the free energy $\mathcal{F}[\rho]$ in the set of densities integrable in $L^1\cap L^m (\Omega)$ sense. General convexity conditions for smooth confinement and interaction potentials allowing for unique minimizers of the free energy $\mathcal{F}[\rho] $ were obtained in \cite{CaMcCVi03}.

The most classical instance corresponds to the choice of attractive Newtonian interaction for $K$ in the whole space without confinement. This case appears both in gravitational collapse and in the mathematical biology literature as a hyperbolic counterpart of the Keller-Segel model, and it is known as the Euler-Poisson system. Conditions on $m$ ensuring the existence of compactly supported H\"older continuous stationary solutions which are minimizers of the free energy $\mathcal{F}[\rho] $ are given in \cite{AB,Str,CS0,CCV,CHVY}. Stability/instability of these stationary solutions together with qualititative properties for the Euler-Poisson system have been analyzed in \cite{M,MP,ELT,LT1,LT2,DLYZ,TW}. A generalization of this case for homogeneous interaction kernels $K(x)=|x|^k/k$, $k\in(-d,0)$, in the whole $d$--dimensional space without confinement has been recently obtained in \cite{CFP,CCH1,CCH2,CHMV}. The balance of forces happen in general for the range $m>1-\tfrac{k}{d}$ regardeless of the mass of the initial data, known as the diffusion-dominated case. The uniqueness of stationary solutions \eqref{balanceforces} has also been established in \cite{CHVY}.

The other interesting case corresponds to quadratic nonlinear diffusion $m=2$ with integrable attractive kernel $K(x)$. This case is directly linked to approximate concentrated repulsion forces by a Dirac kernel in the interaction potential leading to $\varrho(\nabla K\ast \varrho)\simeq \tfrac12 \nabla \varrho^2$. In this case, minimizers of the interaction energy $\mathcal{F}[\rho]$ have been studied in \cite{Bed,BDF,BFH,Kai}. 
In short, they show the existence of a unique compactly supported H\"older continuous stationary solution for all masses if $m>2$ using arguments from \cite{CHVY}. These stationary solutions and their uniqueness, for a given mass and center of mass, are valid both for the whole space and the bounded domain case. The presence of an external potential in the whole space is not needed in order to have a unique stationary state. Adapting \cite{CHVY} to allow for uniformly convex external potentials $\Phi$ in the whole space satisfying the assumptions from the next section is a simple matter.

The energy dissipation \eqref{endis} implies that if we are able to find the suitable analytical framework to show global existence of solutions to the hydrodynamic system \eqref{main_system}, the long-time asymptotics of the system \eqref{main_system} should be governed by the stationary solutions whose densities are given by \eqref{balanceforces}. This is the main goal of this work that we will achieve by finding suitable approximation systems and compactness arguments based on ideas developed by Lions \cite{LI4}, and by Feireisl and collaborators \cite{F_com,EF70,FNP,FNP2002,FP}. Notice that the uniqueness of stationary densities satisfying \eqref{balanceforces} discussed above allows us to uniquely identify long-time asymptotics of hydrodynamic system \eqref{main_system} via compactness methods.

The confining potential $\Phi$ is needed in order confine the mass in the whole space, see \cite{CKT}. Removing $\Phi$ completely is a challenging issue, at least for certain interaction potentials, well known in the aggregation-diffusion literature, see \cite{CCV,CHVY,CS0} for instance. In some cases, quadratic confinement potentials correspond to hydrodynamic systems with expanding self-similar solutions due to repulsive nonlocal interactions in the self-similar variables as in \cite{CCZ}.

Further previous results related to systems of the form \eqref{main_system} include long-time asymptotics and critical thresholds without pressure terms in one dimension \cite{CCTT,CCZ} and weak-strong uniqueness and nonuniqueness results \cite{CFGS} for the corresponding Euler equations. As far as the Navier-Stokes-Poisson (NSP) system is concerned, various questions like existence, long-time behavior, and stability of solutions have been studied during the last 20 years. In the case of Newton potential modelling attractive (gravitational) force in the whole $\R^3$, the existence of global in time weak solutions was proved by Ducomet et al. \cite{DFPS}. Their proof works also in the case of repulsive potential (i.e. Coulomb force) and in this case it has been further studied by Bella in \cite{Bella}. The last of the mentioned results concerns the long-time asymptotics of solutions to NSP. It was shown that these solutions converge to the unique (non-trivial) solution to the stationary problem. In this work we deal with long-time asymptotics for the full hydrodynamic system with nonlocal potential that might be even more singular at the origin. 

The rest of the work is organized as follows. Section \ref{Sec:2} is devoted to the precise hypotheses under which we obtain the main results: existence of solutions for the system  \eqref{main_system} and long-time asymptotics. Section \ref{Sec:3} deals with all the approximation procedures, compactness arguments and estimates needed for the global existence of weak solutions to the system \eqref{main_system}. Further properties related to center of mass and confinement of the mass are obtained for the full space case at the end of this section. Next, in Section \ref{s:long_time}, we analyses the long-time asymptotics showing that the $\omega$-limit set of the constructed weak solutions is determined by the set of stationary solutions satisfying the balance of forces \eqref{balanceforces}. Finally, in Section \ref{Sec:5}, we discuss the extension of our results to the case of system  \eqref{main_system} augmented with terms modelling the alignment.
%%%%%%%%%%%%%%%%%%%%%%%%%%%%%%%%%

\section{Main results: Global Existence and Long-Time Asymptotics}\label{Sec:2}

Let us start by being more precise on the interaction and confinement potentials.
In case $\Omega=\R^3$, we will assume that $\Phi \in W^{1,\infty}_{loc} (\R^3)$ is a  confinement function in the sense that it satisfies that there is $R_o>0$ large enough and $C,\nu>0$ such that
	\begin{equation}\label{ass_Phi_1}
	|\nabla_x \Phi (x) | \leq C\,\Phi(x) \quad \mbox{ for }|x| > R_o, 
	\end{equation}
and
	\begin{equation}\label{ass_Phi_2}
	\lim_{|x| \to \infty } \frac{\Phi(x)}{|x|^{1+ \nu}} = \infty \,.
	\end{equation}
The ensemble of hypotheses $\Phi \in W^{1,\infty}_{loc} (\R^3)$, \eqref{ass_Phi_1} and \eqref{ass_Phi_2} will be referred to as {\bf (HC)}. In case we work in a bounded domain $\Omega$, the confinement potential can be taken as zero if desired. The interaction potential $K(x)$ is assumed to be symmetric and satisfying 
\begin{equation}\label{ass_K}
 K \in L^q(\R^3),\quad  \mbox{ where } \max \left\{\frac{1}{m-1}, 1 \right\}<q < \infty,
 \end{equation}
 and
 \begin{equation}\label{ass_nabla_K}
 \nabla K \in L^q(\R^3),\quad  \mbox{ where } \max \left\{\frac{m}{2(m-1)+\theta},1\right\}<q,  \quad \theta=\min\left\{\frac23m-1,\frac14\right\}.
 \end{equation}
This set of assumptions will be referred as {\bf (HI)}. 

\begin{Remark}\
\begin{itemize}
\item[i)] The hypotheses {\bf (HC)} on the confinement potential are rather natural in order to get a control of the moments of the density by means of the dissipation of free energy \eqref{endis}, see also \cite{CaMcCVi03,CKT}.

\item[ii)] If the interaction potential is smooth $K\in C^2(\R^3)$ with $K, \nabla K, D^2 K\in L^1\cap L^\infty(\R^3)$, the assumptions {\bf (HI)} are met. A particular example is  $K(x)=\pm\exp(-|x|^2)$. By the recent result \cite{Kai}, under the further assumptions that $K$ is radially symmetric $K(x)=k(r)$ with $k$ bounded below, strictly increasing with zero  limit as $r\to \infty$ and $k''(0)<0$, the existence and uniqueness of compactly supported stationary solution to \eqref{balanceforces} with $m=2$ without confinement is ensured. Again, a particular example to consider would be the attractive potential $K(x)=-\exp(-|x|^2)$.

\item[iii)] Assumption \eqref{ass_K} is needed to control the interaction term in the energy \eqref{energytot}. In fact, we obtain that $(K \ast \varrho)\varrho \in L^\infty(0,T; L^1( \R^3))$ by using $\varrho \in L^\infty(0,T; L^1(\R^3) \cap L^m(\R^3))$. More precisely, using the H\"older, Young, and interpolation inequalities respectively we show that
\begin{align}\label{p_ener}
\left|\intO{\vr K\ast\vr}\right| &\leq \|\vr\|_{L^p(\R^3)}\|K\ast\vr\|_{L^{\frac{p}{p-1}}(\R^3)}\leq \|\vr\|^2_{L^p(\R^3)}\|K\|_{L^{q}(\R^3)}\nonumber\\
&\leq C\|\vr\|_{L^1(\R^3)}^{2(1-\alpha)}\|\vr\|_{L^m(\R^3)}^{2\alpha}\|K\|_{L^{q}(\R^3)}\leq C\|\vr\|_{L^m(\R^3)}^{2\alpha}\|K\|_{L^q(\R^3)}\,,
\end{align}
where $q=\frac{p}{2(p-1)}$. Therefore one needs to assume that
\eqh{
1\leq p\leq m,\quad \frac{1}{p}=\frac{1-\alpha}{1}+\frac{\alpha}{m}, \quad 0\leq \alpha\leq 1,\ \ \text{and}\quad 2\alpha<m.}
From this, we deduce two relations between $q$ and $m$, namely $\frac{m}{2(m-1)}\leq q$ and $\frac{1}{m-1}< q$. Therefore, for  $m$ close to  $\frac32$, the second one is most restrictive. Note that for $m<2$, $\frac{1}{m-1}> \frac{m}{2(m - 1)}$, and that for $m\geq 2$, $1\geq \frac{m}{2(m - 1)}$. This leads to our assumption \eqref{ass_K}.

\item[iv)]
If the interaction potential behaves at zero as a power law, in the sense that $K(x)\simeq r^b$ and $|\nabla K(x)|\simeq r^{b-1}$ for small $r$, the integrability hypotheses in {\bf (HI)} are met locally near the origin by taking $q\to 1^+$ as soon as $b>-2$ and $m>\frac{9}5$ after some easy computations using $\frac{m}{2(m-1)+\theta}=1$ in \eqref{ass_nabla_K}. Of course, the potential has to be modified at infinity in order to satisfy the assumptions {\bf (HI)} globally. This case includes a large class of potentials for which the results in \cite{CHVY,CCH1,CCH2} apply giving the existence of radially compactly supported stationary solutions to \eqref{balanceforces} without confinement.  This class therefore allows for potentials that are locally even more singular than Newtonian interaction. Our present results do not include the case of the purely attractive/repulsive Newtonian interaction for which uniqueness of the solution to \eqref{balanceforces} is known \cite{CCV,CHVY},  see also \cite{Bella}.

\item[v)] Notice that by taking in \eqref{ass_nabla_K} $m\to\frac32^+$, $\theta\to 0^+$, and therefore $q> \frac 32$. 
\end{itemize}
\end{Remark}

We will denote by $\D$ the space of smooth compactly supported functions, by 
 $D^{1,2}(\Omega)$ -- the space of locally $L^1$-integrable functions with gradient in $L^2(\Omega)$.
 The symbol $C_{weak}([0,T];L^q(\Omega))$ stands for the space of all  vector-valued functions on $[0,T]$ ranging in $L^q(\Omega)$ continuous with respect to the weak topology.
Let us first define the concept of weak solution that we deal with in this work.

\begin{Definition}\label{defsol} Let $\Omega =\R^3$ or be a bounded smooth set in $\R^3$. 
Given the density $\varrho \in L^\infty(0,T; L^1_+ \cap L^m(\Omega))$ with  $\varrho \in C_{weak}([0,T]; L^m(\Omega))$ and the velocity field $\varrho\vu \in C_{weak}([0,T]; L^{\frac{2m}{m+1}}(\Omega))$ with  $\vu \in L^2(0,T; W^{1,2}_0(\Omega))$ in the case of bounded domain or $\vu \in L^2(0,T; D^{1,2}(\Omega))$ in the case $\Omega = \R^3$, we say that the pair $(\varrho,\vu)$ is a bounded free energy weak solution to the hydrodynamic system \eqref{main_system} if 
\begin{itemize}
\item[i)] $\varrho \geq 0$ is a renormalized solution of continuity equation \eqref{main_system}$_1$ on $(0,T) \times \Omega$: 
for any test function $\varphi \in \D ( [0,T) \times \Omega)$, any $T>0$ and any $b \in C^1(\R)$ s.t. $b'(z) =0$ for all $z$ large enough, say $z > M_b$, the following holds
	\begin{equation}\label{weak_1}
	\int_0^T \int_{\Omega} \left( b(\varrho) \partial_t \varphi + (b(\varrho)  \vu) \cdot \nabla_x \varphi 
	- (b'(\varrho)\varrho  - b(\varrho) ) \Div \vu \varphi  \right)
	\dxdt = - \int_{\Omega} b(\varrho_0) \varphi(0,\cdot) \dx .
	\end{equation} 
 Moreover, \eqref{weak_1} holds on $\R^3$ provided $\varrho$, $\vu$ is extended to be zero on $\R^3 \setminus \Omega$ in case of $\Omega$ bounded.

\item[ii)] The equation \eqref{main_system}$_2$ is satisfied in distributional sense, i.e.
	\begin{align*}%\label{weak_2}
	\int_0^T\!\!\!\int_{\Omega} \left( \varrho \vu \cdot \partial_t\vvarphi  + \varrho \vu \otimes \vu : \nabla_x \vvarphi  
	+ a \varrho^m \Div \vvarphi\right) \dxdt =  \nonumber & 
	\int_0^T\!\!\!\int_{\Omega} \big( \mu \nabla_x \vu : \nabla_x \vvarphi   + (\lambda + \mu) \Div \vu \Div \vvarphi \big) \dxdt 
	\\ \nonumber & 
	+  \int_0^T\!\!\!\int_{\Omega} \left( ( \nabla_x K \ast \varrho + \nabla_x \Phi)  \varrho + \varrho\vu \right) \cdot \vvarphi 
	\dxdt
	\\ & 
	- \int_{\Omega} \vm_0 \cdot \vvarphi(0) \dx,
	\end{align*}
for any test function $\vvarphi \in \D ([0,T) \times \Omega ;\R^3)$ and any time $T>0$.

\item[iii)] The following energy inequality
$$
E(t) + \int_0^t \int_{\Omega}\left[ \mu |\nabla \vu |^2 + (\lambda+ \mu) | \Div \vu|^2 +   \varrho |\vu|^2  \right]
\dx \dt \leq E(0)
$$
holds for a.e. $t\in (0,T)$,  where
 	$$
 	E (t) = \int_{\Omega} \left( \frac{1}{2} \varrho(t) |\vu (t) |^2 + \frac{a}{m-1} \varrho^m(t)  + \frac{1}{2} (K \ast \varrho) \varrho +  \varrho \Phi  \right) \dx .
 	$$
\end{itemize}
\end{Definition}

Our main result concerning global existence of weak solutions in the whole space case reads as follows. 
	
\begin{Theorem}\label{t_existence}
Let $m>\frac{3}{2}$ and $\Omega=\R^3$. Assume that the interaction potential $K$ and the confinement potential $\Phi$ satisfy {\bf (HI)} and {\bf (HC)} respectively, and that the initial data satisfy \eqref{initial_data}-\eqref{com_con} together with $\varrho_0 \Phi  \in L^1(\R^3)$.  

Then there exists a global in time bounded free energy weak solution to the system \eqref{main_system}, \eqref{initial_data}, in the sense of Definition  {\rm\ref{defsol}}. Moreover, for $\theta = \min \{ \frac{2}{3} m - 1, \frac{1}{4} \}$, $\varrho$ satisfies additionally
	\begin{equation}\label{t1_3_bis}
	 \varrho \in L_{loc}^{m + \theta} ((0,T) \times \R^3).
	\end{equation}
\end{Theorem}

\begin{Remark}\label{Rem_2.2}
Assumption \eqref{ass_nabla_K} will be used to control the interaction term for obtaining the higher integrability of the density \eqref{t1_3_bis} locally in space, and the weak sequential stability of the nonlocal term, see Lemma \ref{lem_theta} below.
\end{Remark}

In the bounded domain case, there is no need for the confinement potential, the assumptions \eqref{ass_Phi_1},  \eqref{ass_Phi_2} do not have to be satisfied, in particular, we may take $\Phi \equiv 0$. Moreover, for the bounded domain the hypotheses {\bf (HI)} on the nonlocal force $\nabla K\ast \varrho$ can also be reduced to their local versions. We will refer to that set of hypotheses as {\bf (HI)$_{loc}$}. This is due to the fact that for densities supported inside a ball of radius $R$ we only need values of $K(x)$ or their derivatives inside a ball of radius at most $2R$. Therefore, given a potential $K$ satisfying {\bf (HI)$_{loc}$}, we can change $K$ outside of the ball of radius $2R$ to satisfy {\bf (HI)}. Therefore without loss of generality, we can assume in the theorem {\bf (HI)$_{loc}$} although we will use in the proofs {\bf (HI)} if needed. Then the main result for the bounded domain is:

\begin{Theorem}\label{t_existenceb}
Let $m>\frac{3}{2}$ and $\Omega$ a smooth bounded domain in $\R^3$. Assume the interaction potential $K$ satisfies {\bf (HI)$_{loc}$}, $\Phi \in W^{1,\infty} (\Omega)$, and that the initial data satisfy \eqref{initial_data}-\eqref{com_con}. 

Then there exists a global in time bounded free energy weak solution to the system \eqref{main_system},  \eqref{bc_intro}, \eqref{initial_data}, in the sense of Definition {\rm\ref{defsol}}. Moreover, for $\theta = \min \{ \frac{2}{3} m - 1, \frac{1}{4} \}$, $\varrho$ satisfies additionally
	\begin{equation*}%\label{t1_3_bis_b}
	 \varrho \in L^{m + \theta} ((0,T) \times \Omega).
	\end{equation*}
\end{Theorem}

Now, we can also discuss the main result concerning the long-time asymptotics of the hydrodynamic system \eqref{main_system}. Given a curve of weak solutions in $\varrho\in C_{weak}([0,\infty); L^m (\Omega))\cap L^\infty((0,\infty);L^1\cap L^m (\Omega))$, we define its $\omega$-limit set $\omega(\varrho)$ in $L^m(\Omega)$ as the set of all possible accumulation points of the curve as $t\to \infty$, i.e.,
$$
\omega(\varrho) = \left\{ \bar\varrho\in L^m(\Omega) \mbox{ such that } \exists \{t_n\}_{n\in\N}\nearrow \infty : \{\varrho(t_n,x)\}_{n\in\N}\to \bar\varrho(x) \mbox{ strongly in } L^m(\Omega) \right\}\,.
$$

\begin{Theorem}\label{t_longtime}
The $\omega$-limit set $\omega(\varrho)$ associated to global weak solutions $(\varrho(t),\vu(t))$ to the hydrodynamic system \eqref{main_system} obtained in Theorems {\rm\ref{t_existence}} and {\rm\ref{t_existenceb}} consists of stationary solutions with zero  momentum and densities with the same initial mass $M_0=\|\varrho_0\|_{L^1(\Omega)}$, satisfying the balance of forces relation \eqref{balanceforces} in $\D'(\Omega)$.
%%%% This contains the no flux boundary condition implicitly for $\Omega$  bounded
Moreover, we have
$$\lim_{t\to\infty} \|\Grad\vu(t)\|_{L^2 (\Omega)}=0 ,$$
and
\begin{equation}\label{u_to_zero_a}
\lim_{t\to\infty} \|\vu(t)\|_{L^6 (\Omega)}=0\, \mbox{ in case of bounded domain }\Omega
\end{equation}
\begin{equation}\label{u_to_zero_b}
\lim_{t\to\infty} \int_\Omega \varrho (t) | \vu(t) |^2 \dx =0 \quad \mbox{ for } \Omega = \R^3 \mbox{ and } \Omega \mbox{ bounded}.
\end{equation}
Additionally, if the solution to \eqref{balanceforces} with mass $M_0$ and  zero-center of mass density is uniquely given by $\varrho_s$, and there exists $c_0\in\R^3$, such that
$$
\lim_{t\to\infty} \int_{\R^3} x\varrho(t)\dx = c_0\,,
$$
then
$$
\lim_{t\to\infty} \|\varrho(t,x)-\varrho_s(x-c_0)\|_{L^m(\R^3)}=0\,.
$$
\end{Theorem}

\begin{Remark}
According to \cite{Kai}, if $m>2$ and $K(x)=1- \exp({- \frac{|x|^2}{2}})$ the solution to \eqref{balanceforces} is unique up to translation and given by a H\"older continuous compactly supported radially decreasing profile $\varrho_s$. Therefore, the convergence towards the unique steady state holds for the corresponding hydrodynamic system \eqref{main_system}. For a more general set of assumptions on $K$ and $m$ under which the uniqueness up to translations of solutions to \eqref{balanceforces} holds, we refer to \cite{Kai}.
\end{Remark}

%%%%%%%%%%%%%%%%%%%%%%%%%%%%%%%%%%%%%%%%

\section{Proof of global existence}\label{Sec:3}

In the following section we prove existence of solutions to the system \eqref{main_system}. Our proof is based on the existence of solutions from \cite{FNP}. The result therein holds for the Navier-Stokes system on bounded domain and without damping, confinement and nonlocal terms. In this section, we concentrate on the whole space case, and therefore we introduce three levels of approximations:
	\begin{itemize}
	\item  We introduce an approximation by bounded domains: balls of radius $r$, that will converge later to whole space $\R^3$ as $r \to \infty$.
	\item As in \cite{FNP}, we introduce artificial viscosity term in the continuity equation related to the parameter $\ep$, which later converges to zero.
	\item As in \cite{FNP}, we introduce the artificial pressure term related to the parameter $\delta$, which later converges to zero.
	\end{itemize}
In order to prove the existence on bounded domains we can skip the first point of the approximation and it will be easy to observe that confinement is not needed, since \eqref{ass_Phi_1}, \eqref{ass_Phi_2} is crucial only for the case of unbounded domain, in particular  as $r\to \infty$.

\subsection{Approximation}

In order to prove existence of global in time weak solutions to \eqref{main_system} let us start with the following approximation: 
%	\begin{equation}\label{Omega_r}
%	\begin{split}    
%	& \mbox{ Let }\Omega_r \mbox{  be bounded domain in } \R^3 \mbox{ with } C^{2+ \nu}\ (\nu >0) \mbox{ boundary for each }r>0.
%	\\& 
%	 \mbox{ For the family } \{ \Omega_r\}_{r\geq1} \mbox{ we assume that } \Omega_{r-1} \subset \Omega_r \mbox{ and }\bigcup_{r=1}^\infty \Omega_r = \R^3.
%	\end{split}
%	\end{equation}
	\begin{equation}\label{Omega_r}
	\begin{split}    
	 \mbox{ Let }\Omega_r=\{x\in \R^3: \ |x|<r\}, \quad r>0. 
	\end{split}
	\end{equation}
	
	 For $\delta>0$, $\ep>0$ and $\beta >\max\{4,m\}$ we introduce the following system which approximates \eqref{main_system}
	\begin{align}
	\partial_t \varrho + \Div (\varrho \vu) 
	 = & \ep \Delta \varrho \quad \mbox{ in } (0,T) \times \Omega_r \,,\nonumber
	\\
	\partial_t(\varrho \vu) + \Div (\varrho \vu \otimes \vu ) + \nabla_x (a \varrho^m + \delta \varrho^\beta) + \ep \nabla_x \varrho \nabla_x \vu 
	 =\,&  \mu \Delta \vu + (\lambda + \mu) \nabla_x (\Div \vu)  \label{usmiech_a1}\\
	 & - \varrho\vu - (\nabla_x K \ast \varrho  + \nabla_x \Phi ) \varrho \quad \mbox{ in }  (0,T) \times \Omega_r\nonumber
	\end{align}
supplemented with  the  boundary conditions on $\partial \Omega_r$
	\begin{equation}\label{bc_a1}
	\nabla_x \varrho \cdot \vn |_{\partial \Omega_r} = 0, \quad \quad \vu |_{\partial \Omega_r} = 0,
	\end{equation}
where $\vn$ is an outer to $\Omega_r$ normal unit vector. Parameters $\mu$, $\lambda$ and $m$ satisfy \eqref{lambda_mu_m}. Moreover we complete the system with the initial data:
	\begin{equation}\label{id_a1}
	\begin{split}
	& \varrho_\delta(0) = \varrho_{0,\delta} \in C^{2+\nu}(\overline{\Omega_r}), \quad 0 < \underline{\varrho} \leq \varrho_{0,\delta} (x) \leq \overline{\varrho} < \infty \mbox{ in }\Omega_r, 
	\\ &
	(\varrho_\delta \vu_\delta)(0) = \vm_{0,\delta} \in C^2(\overline{\Omega_r}).
	\end{split}
	\end{equation}

\subsubsection{Existence of solutions to the approximation}
The proof of  existence of solutions to the system \eqref{usmiech_a1}, \eqref{bc_a1}, \eqref{id_a1} with $\Omega_r$ satisfying \eqref{Omega_r}
follows the blueprint of \cite{FNP}. We proceed analogously, since the new terms (damping, confinement, and nonlocal) do not cause any additional difficulties as proved below. Namely, one can construct Faedo-Galerkin approximation scheme in the following way (see \cite[Section~2.2]{FNP} for details):
\begin{enumerate}
\item Velocity field is spanned in finitely dimensional space of functions of suitable regularity (eigenfunctions of the Laplacian).
\item The continuity equation \eqref{usmiech_a1}$_1$ with initial data \eqref{id_a1} has then classical solutions.
\item Let us call $k$-approximation any such constructed approximation. In order to use fixed point argument which provides existence of $k$-approximations, we notice that
due to assumptions on the initial data, on $K$, and $\Phi$,  the nonlocal and confinement terms are bounded.
\item  In order to prove weak sequential stability, we  discuss only  the novelty related to the nonlocal term. 
The analysis of $k-$approximation provides that $\varrho_k \to \varrho$ strongly in $L^4((0,T)\times \Omega_r)$ and a.e. in $(0,T)\times \Omega_r$. Since $\nabla_x K \in L^q(\R^3)$ with $q\geq 1$, prolonging $\varrho_k$ by zero outside of $\Omega_r$, $\nabla_x K \ast \varrho_k$ is bounded in $L^4((0,T)\times \R^3)$, $(\nabla_x K \ast \varrho_k) \varrho_k $ is bounded in $L^2((0,T)\times\Omega_r)$. To identify the   weak limit of $(\nabla_x K \ast \varrho_k) \varrho_k $ let us notice that
for any $\vc\phi \in L^2 ((0,T) \times \Omega_r;\R^3)$
\eqh{
	&\int_0^T\int_{\Omega_r}
	(\nabla_x K \ast \varrho_k) \varrho_k \cdot \vc\phi - (\nabla_x K \ast \varrho ) \varrho \cdot\vc\phi \dxdt\\
	&\quad= \int_0^T\int_{\Omega_r}
	(\nabla_x K \ast ( \varrho_k - \varrho)) \varrho_k \cdot \vc\phi - (\nabla_x K \ast \varrho ) 
	(\varrho_k - \varrho) \cdot \vc\phi \dxdt .
}
	By the above properties one easily concluded that both terms on the RHS vanish as $k \to \infty$, consequently $(\nabla_x K \ast \varrho_k) \varrho_k$ converges weakly to $(\nabla_x K \ast \varrho) \varrho$ in $L^2((0,T)\times \Omega_r)$.
\item The confinement term does not bring new difficulties since $\Phi \in W^{1,\infty}_{loc}(\R^3)$, similarly as the dumping term.
\end{enumerate}
 
\subsubsection{Vanishing artificial viscosity limit }\label{sec_vanish_epsilon}

Next goal is to pass to the limit in \eqref{usmiech_a1} with $\ep \to 0$. Let us denote weak solutions to \eqref{usmiech_a1}, \eqref{bc_a1}, \eqref{id_a1} by $\{\varrho_\ep, \vu_\ep\}_{\ep>0}$. Then the following proposition holds:

\begin{Proposition}\label{prop_2.1}
Suppose $\beta > \max \{ 4, m\}$. Let $\Omega_r \subset \R^3$ satisfy \eqref{Omega_r}, $\varrho_{0,\delta}$ and $\vm_{0,\delta}$ satisfy \eqref{bc_a1}, \eqref{id_a1}.  Let $K$ satisfy  {\bf (HI)}, $\Phi \in W^{1,\infty} (\Omega_r)$.  Let $\ep >0$, $\delta>0$.
Then there exists a weak solution $\{ \varrho_\ep, \vu_\ep \}_{\ep >0}$ of the problem \eqref{usmiech_a1}, \eqref{bc_a1}, \eqref{id_a1} such that
	\begin{equation}\label{p21_2}
	\sup_{t\in[0,T]} \| \varrho_\ep (t) \|^m_{L^m(\Omega_r)} \leq c(E_{\delta,r,\ep}(0),  K, \Phi)\,,
	\end{equation}
	\begin{equation}\label{p21_31}
	\delta \sup_{t\in[0,T]} \| \varrho_\ep (t) \|^\beta_{L^\beta(\Omega_r)} 
	\leq c (E_{\delta,r,\ep}(0),  K, \Phi)\,,
	\end{equation}
	\begin{equation}\label{p21_3}
	\sup_{t\in[0,T]} \| \sqrt{\varrho_\ep (t)} \vu_\ep   \|^2_{L^2(\Omega_r)} 
	\leq c (E_{\delta,r,\ep}(0),  K, \Phi)\,,
	\end{equation}
	\begin{equation}\label{p21_4}
	\int_0^T  \|  \vu_\ep   \|^2_{L^2(\Omega_r)} +  \|  \nabla_x \vu_\ep   \|^2_{L^2(\Omega_r)} \dt
	\leq c (E_{\delta,r,\ep}(0),  K, \Phi)\,,
	\end{equation}
	\begin{equation}\label{p21_4_2}
	\sup_{t\in[0,T]} \| \varrho_\ep (t) \Phi  \|_{L^1(\Omega_r)} 
	\leq c (E_{\delta,r,\ep}(0),  K, \Phi)\,,
	\end{equation}
	\begin{equation}\label{p21_5}
	\ep \int_0^T  \| \nabla_x \varrho_\ep   \|^2_{L^2(\Omega_r)} \dt 
	\leq c(\beta,\delta,\varrho_0,\vm_0,K)\,,
	\end{equation}
	\begin{equation*}\label{p21_6}
	\begin{split}
	\frac{\rm d}{{\rm d}t} \int_{\Omega_r} &
	\left[ \frac{1}{2} \varrho_\ep | \vu_\ep |^2 + \frac{a}{m-1} \varrho_\ep^m  
	+ \frac{\delta}{\beta -1}\varrho^{\beta}_\ep   \right] \dx \\
	&+ \int_{\Omega_r} \ep \left[ a\gamma |\nabla_x \varrho_\ep |^2 \varrho^{\gamma-2} + \delta \beta |\nabla_x \varrho|^2 \varrho^{\beta-2} \right] \dx
	+
	\int_{\Omega_r}\left[\mu | \nabla_x \vu_\ep |^2 + (\lambda + \mu) | \Div \vu_\ep |^2 
	+ \varrho_\ep |\vu_\ep|^2\right] \dx\\
	&
	\leq -  \int_{\Omega_r}\left[ \varrho_\ep \nabla_x \Phi \cdot  \vu_\ep + \varrho_\ep (\nabla_x K \ast \varrho_\ep) \cdot  \vu_\ep \right] \dx
	\end{split}
	\end{equation*}
holds in $\D'(0,T)$,  alongside with its integrated version
\begin{equation*}\label{p21_6b}
	\begin{split}
	E_{\delta,r,\ep}(t)
	&+ \int_0^t\int_{\Omega_r} \ep \left[ a\gamma |\nabla_x \varrho_\ep |^2 \varrho^{\gamma-2} + \delta \beta |\nabla_x \varrho|^2 \varrho^{\beta-2} \right] \dx\dt\\
	&+
	\int_0^t\int_{\Omega_r}\left[\mu | \nabla_x \vu_\ep |^2 + (\lambda + \mu) | \Div \vu_\ep |^2 
	+ \varrho_\ep |\vu_\ep|^2\right] \dx\dt\\
	&
	\leq E_{\delta,r,\ep}(0)-  \int_{\Omega_r}\left[ \varrho_\ep \nabla_x \Phi \cdot  \vu_\ep + \varrho_\ep (\nabla_x K \ast \varrho_\ep) \cdot  \vu_\ep \right] \dx
	\end{split}
	\end{equation*}
 for a.a. $t\in(0,T)$, where
	\begin{equation*}
	E_{\delta,r,\ep}(t)= 
	\int_{\Omega_r } \left[ \frac{1}{2} \varrho_\ep | \vu_\ep |^2 + \frac{a}{m-1} \varrho_\ep^m  
	+ \frac{\delta}{\beta -1}\varrho^{\beta}_\ep +  \right] \dx,
	\end{equation*}
and
	\begin{equation*}
	E_{\delta,r,\ep}(0)= 
	\int_{\Omega_r } \left[  
	\frac{1}{2} \frac{|\vm_{0,\delta}|^2}{\varrho_{0,\delta}}  + \frac{a}{m-1} \varrho_{0,\delta}^{m} 
	+ \frac{\delta}{\beta-1}\varrho_{0,\delta}^\beta
	% + \frac{1}{2}(K \ast \varrho_{0,\delta})\varrho_{0,\delta} + \varrho_{0,\delta} \Phi 
	 \right] \dx.
	\end{equation*}
There exists also $s>1$ s.t. $\partial_t \varrho_\ep$, $\Delta \varrho_\ep \in L^s((0,T)\times \Omega_r)$ and continuity equations is satisfied a.e. in $(0,T)\times \Omega_r$, and
		\begin{equation}\label{p21_8}
		\varrho_\ep \to \varrho \quad \mbox{ in }C_{weak}((0,T); L^\beta (\Omega_r)).
		\end{equation}
	Moreover there exists a constant $c$ independent of $\ep$ s.t.
		\begin{equation}\label{p21_7}
		\| \varrho_\ep \|_{L^{m+1}((0,T)\times \Omega_r)} + \| \varrho_\ep \|_{L^{\beta+1}((0,T)\times \Omega_r)}
		\leq c(\delta,\varrho_{0,\delta},\vm_{0,\delta},K) .
		\end{equation}
\end{Proposition}
Details of  the proof can be found in \cite[Section~2]{FNP}. One needs only to modify it  due to presence of damping, nonlocal, and confinement terms. 
Estimates \eqref{p21_2} -- \eqref{p21_4} are a consequence of the energy inequality, application of the Poincar\'e and the Young inequality, and assumptions on $K$ and $\Phi$.  
To justify the estimate \eqref{p21_7} an additional step is to give uniform bounds for damping, nonlocal, and confinement terms in the momentum equation tested by $\psi(t) \B[\varrho_\ep - k_0]$, where $\B$ denotes Bogovski operator, $k_0 = \frac{1}{|\Omega_r|} \int_{\Omega_r} \varrho_\ep (x)\dx$, $\psi \in \D(0,T)$, $0\leq \psi \leq 1$. Then as $\B[\varrho_\ep - k_0]$ is bounded in $L^\infty(0,T;L^\infty(\Omega_r))$ as $\beta > 4 $, \eqref{p21_2} - \eqref{p21_4_2} hold, and by assumptions on $\nabla_x K$  and on $\Phi$ we find that
 	\begin{equation*}
	\left|\int_0^T \psi \int_{\Omega_r} \varrho_\ep \vu_\ep \B[\varrho_\ep - k_0] \dxdt\right|
	\leq  c_1 \int_0^T \| \sqrt{\varrho_\ep}\|_{L^2(\Omega_r)} \| \sqrt{\varrho_\ep} \vu_\ep\|_{L^2(\Omega_r)}
	\| \B[\varrho_\ep - k_0] \|_{L^\infty(\Omega_r)} \dt \leq c_2,
	\end{equation*}
 	\begin{equation*}
	\left|\int_0^T \psi \int_{\Omega_r} (\nabla_x K \ast \varrho_\ep) \varrho_\ep \B[\varrho_\ep - k_0] \dxdt \right|
	\leq  c_1 \int_0^T \| \nabla_x K\ast \varrho_\ep\|_{L^2(\Omega_r)} \| \varrho_\ep\|_{L^2(\Omega_r)}
	\| \B[\varrho_\ep - k_0] \|_{L^\infty(\Omega_r)} \dt \leq c_2,
	\end{equation*}
	\begin{equation*}
	\begin{split}
	& \left| \int_0^T \psi \int_{\Omega_r} (\nabla_x \Phi) \varrho_\ep \B[\varrho_\ep - k_0] \dxdt \right|
	 \\
	& \qquad \qquad \leq  c_1 \int_0^T \left(\| \nabla_x \Phi \|_{L^{m'}(\Omega_r \cap B_{R_o})} \| \varrho_\ep\|_{L^m(\Omega_r \cap B_{R_o})} + \| \Phi \varrho_\ep\|_{L^1(\Omega_r  \setminus B_{R_o})}\right)
	\| \B[\varrho_\ep - k_0] \|_{L^\infty(\Omega_r)} \dt \leq c_2
	\end{split}
	\end{equation*}
where $B_{R_o}$ centered at zero with radius ${R_o}$ given in \eqref{ass_Phi_1}.	

Let us notice that as $\varrho_\ep$ satisfies continuity equation in a weak sense
then  we get that for all $\phi \in C^\infty_c (\Omega_r)$, the family 
$\left\{ \int_{\Omega_r} \varrho_\ep \phi \dx \right\}_{\ep>0} (t) $
forms a bounded and equicontinuous sequence in $C[0,T]$ (since by \eqref{p21_2}, \eqref{p21_3}, $\varrho_\ep$ bounded in $L^\infty(0,T;L^m(\Omega_r))$ and $\varrho_\ep \vu_\ep$ bounded in $L^\infty(0,T;L^{\frac{2m}{m+1}}(\Omega_r))$). Then, by the Arzela-Ascoli theorem 
$$
\int_{\Omega_r} \varrho_\ep \phi\dx \to \int_{\Omega_r} \varrho_\ep \phi \dx \mbox{ in }C([0,T]) \mbox{ for any } \phi \in C^\infty_c (\Omega_r).
$$ 
As \eqref{p21_2} holds,  the above convergence extends to each $\phi \in L^{m'}(\Omega_r)$ (see \cite[Corollary~2.1]{EF70}) and therefore we get
\begin{equation}\label{rho_ep_weak_CLm}
	\varrho_\ep \to \varrho \quad \mbox{ in }C_{weak}([0,T];L^m(\Omega_r)).
\end{equation}
So, by \eqref{p21_31} we get that \eqref{p21_8} holds as well. 
Similar arguments as for \eqref{rho_ep_weak_CLm} applied to $\varrho_\ep \vu_\ep$ and momentum equation  combined with  \eqref{p21_3}, \eqref{p21_4} provide 
	\begin{equation}\label{rhou_lim_ep}
	\varrho_\ep \vu_\ep \to \varrho \vu \quad \mbox{ in } C_{weak}([0,T]; L^{\frac{2m}{m+1}}(\Omega_r)).
	\end{equation}
Combining this with the weak convergence of the gradient of $\vu_\ep$ allows us to pass to the limit in the convective term, using the Div-Curl argument.

With the Proposition~\ref{prop_2.1} at hand, we may pass with $\ep \to 0$ in the continuity and momentum equations \eqref{usmiech_a1} and obtain (in a weak sense):
	\begin{align}
	\partial_t \varrho + \Div (\varrho \vu) 
	& = 0 \quad \mbox{ in } (0,T) \times \Omega_r\,,\label{usmiech_a2}
	\\
	\partial_t(\varrho \vu) + \Div (\varrho \vu \otimes \vu ) + \nabla_x \overline{p}
	& = \mu \Delta \vu + (\lambda + \mu) \nabla_x (\Div \vu)  - \varrho\vu - \overline{(\nabla_x K \ast \varrho)\varrho}  - \varrho \nabla_x \Phi   \mbox{ in }  (0,T) \times \Omega_r. \nonumber
	\end{align}
Here, $\overline{p}$ is a weak limit of $a \varrho^m_\ep + \delta \varrho^\beta_\ep$, namely by 
\eqref{p21_7}
	\eq{\label{lim_den}
	a \varrho^m_\ep \weak \overline{a \varrho^m} \mbox{ weakly in }L^{\frac{m+1}{m}}((0,T) \times \Omega_r), \\
	 \delta\varrho^\beta_\ep \weak \overline{\delta\varrho^\beta} \mbox{ weakly in }L^{\frac{\beta+1}{\beta}}((0,T) \times \Omega_r) .
	}
and so $ \overline{p} = \overline{a \varrho^m} + \overline{\delta\varrho^\beta} $. By
$\overline{(\nabla_x K \ast \varrho)\varrho}$ we denote the weak limit of  $(\nabla_x K \ast \varrho_\ep)\varrho_\ep$ in $L^s((0,T) \times \Omega_r)$ for some $s>1$, due to \eqref{p21_2} and assumptions {\bf (HI)} on $K$. 

To show that
	\begin{equation}\label{p_lim_ep}
	 \overline{p} = a \varrho^m + \delta \varrho^\beta \,,
	 \end{equation}
we will follow arguments used in \cite[Section~3.4]{FNP}, in particular let us state the following result concerning weak sequential stability of the so-called effective viscous flux:
	\begin{Lemma}\label{lem_3.2}
	Let $\{ \varrho_\ep, \vu_\ep\}_{\ep>0}$ be a weak solutions to \eqref{usmiech_a1}, \eqref{bc_a1}, \eqref{id_a1} as in Proposition~\ref{prop_2.1} and $(\varrho,\vu)$ be a weak solution to \eqref{usmiech_a2}. Then
	$$
	\lim_{\ep\to 0^+} \int_0^T \psi \int_{\Omega_r} \phi ( a \varrho^m_\ep + \delta \varrho^\beta_\ep - (\lambda + 2 \mu) \Div \vu_\ep)
	 \varrho_\ep \dxdt
	 = \int_0^T\psi \int_{\Omega_r} \phi ( \overline{p} - (\lambda + 2 \mu) \Div \vu)
	 \varrho \dxdt
	$$
for all $\psi \in \D(0,T)$, $\phi \in \D(\Omega_r)$.
	\end{Lemma}
The proof of Lemma~\ref{lem_3.2} is almost the same as for \cite[Lemma~3.2]{FNP}. We only need to take care of damping, nonlocal and confinement terms when testing the momentum equation by properly chosen functions, passing to the limit and checking if these terms converge properly to their counterparts.  To this end let us denote
	$$\varphi = \psi(t) \phi(x) \A [\varrho_\ep]$$
  where  $\psi \in \D(0,T)$, $\phi \in \D(\Omega_r)$, $\varrho_\ep$ is extended by zero outside of $\Omega_r$ and the inverse divergence operator $\A$ is defined by
  	\begin{equation}\label{inv-div}
	\A[v] = \{\A_i\}_{i=1,2,3},\  \A_i= \Delta^{-1} (\partial_{x_i} v)\,.
	\end{equation}           
In particular, for ${\mathcal{A}}_i [h]  = \Delta^{-1} \partial_j h$ by Lizorkin theorem we have
	\begin{equation}\label{Lizorkin}
 	{\mathcal{A}}_i : L^p( \R^3)  \to L^q (\R^3) \mbox{
 	with } q= \frac{3p}{3 - p },\ 1 < p \leq q < \infty,\ p<3.
	 \end{equation}
	
Notice that $\varphi$ defined above is a proper test function for momentum equation, particularly for nonlocal, damping, and confinement  terms. Namely $\A [\varrho_\ep] \in L^\infty(0,T;L^\infty(\R^3))$, as  $\A [\varrho_\ep] \in L^\infty(0,T;W^{1,\beta}(\R^3))$ and $\beta>4$.
The properties of operator $\A$ give 
	$$\| \A [v] \|_{W^{1,s}(\Omega_r)} \leq c(s,\Omega_r) \| v \|_{L^s(\R^3)} \quad \mbox{ with } 1<s<\infty.$$
Then by \eqref{p21_8} and \eqref{p21_7} we get
	\begin{equation}\label{A_lim_ep}
	\A [\varrho_\ep] \to \A[\varrho] \quad \mbox{ in } C(\overline{(0,T)\times \Omega_r} ) .
	\end{equation}
Hence by \eqref{A_lim_ep} and \eqref{rhou_lim_ep}, we deduce
	$$
	\int_0^T \int_{\Omega_r} \varrho_\ep \vu_\ep \cdot  \psi(t) \phi(x) \A [\varrho_\ep] \dxdt \to \int_0^T \int_{\Omega_r} \varrho \vu \cdot \psi(t) \phi(x) \A [\varrho] \dxdt \quad \mbox{ as } \ep\to 0.
	$$
By \eqref{A_lim_ep}, one obtains 
		$$
	\int_0^T \int_{\Omega_r} (\nabla_x K \ast \varrho_\ep)\varrho_\ep \cdot \psi(t) \phi(x) \A [\varrho_\ep] \dxdt \to \int_0^T \int_{\Omega_r} \Ov{(\nabla_x K \ast \varrho)\varrho}  \cdot \psi(t) \phi(x) \A [\varrho] \dxdt \quad \mbox{ as } \ep\to 0.
	$$
Note that this term is precisely what we obtain by testing the limit momentum equation by the function $ \psi(t) \phi(x) \A [\varrho]$.
Next, since $\Phi \in W^{1,\infty}_{loc} (\R^3)$,  \eqref{A_lim_ep}  and \eqref{lim_den} provide that 
		$$
	\int_0^T \int_{\Omega_r} \varrho_\ep \nabla_x \Phi \cdot   \psi(t) \phi(x) \A [\varrho_\ep] \dxdt \to \int_0^T \int_{\Omega_r} \varrho   \nabla_x \Phi \cdot \psi(t) \phi(x) \A [\varrho] \dxdt \quad \mbox{ as } \ep\to 0.
	$$
\bigskip

Lemma~\ref{lem_3.2} is crucial to provide strong convergence of the density sequence.  For the rest of the details we refer to  \cite{FNP}. Just shortly, it is based on the fact that $P(z) = a z^m + \delta z^\beta$ is a monotone function, and one can use a Minty type arguments to prove a.e. convergence of the density sequence
	$$
	\varrho_\ep \to \varrho \quad \mbox{ a.e. in } (0,T)\times \Omega_r,
	$$
and that \eqref{p_lim_ep} is satisfied. Moreover, thanks to uniform estimates on the density, we have that 
$$
	\varrho_\ep \to \varrho \quad \mbox{ strongly in } L^p((0,T)\times \Omega_r),\quad 1\leq p<\beta+1.
	$$
	So, using  {\bf (HI)} we identify also
	$$ 
	\Ov{(\nabla_x K \ast \varrho)\varrho} = (\nabla_x K \ast \varrho)\varrho.
	$$
In this way we are able to conclude the following result.

\begin{Proposition}\label{prop31}
Let $\Omega_r \subset \R^3$ be bounded domain of the class  $C^{2+ \nu}$, $\beta > 4.$  Let $K$ satisfy  {\bf (HI)} and  $\Phi$ satisfy ${\bf (HC)}$.
Let \eqref{id_a1} be satisfied. Then there exists a free energy weak solution $( \varrho, \vu )$ to the problem
	 	\begin{align}
	\partial_t \varrho + \Div (\varrho \vu) 
	 = \,& 0 \quad \mbox{ in } (0,T) \times \Omega_r, \nonumber
	\\
	\partial_t(\varrho \vu) + \Div (\varrho \vu \otimes \vu ) + \nabla_x (a \varrho^m + \delta \varrho^\beta) 
	 = \, & \mu \Delta \vu + (\lambda + \mu) \nabla_x (\Div \vu) \label{usmiech_a3}
	 \\ & - \varrho\vu - (\nabla_x K \ast \varrho + \nabla_x \Phi)\varrho   \quad \mbox{ in }  (0,T) \times \Omega_r, \nonumber
	\\ 
	\vu |_{\partial\Omega_r}  = & \,0 \nonumber
	\end{align}
for any fixed $r >1$ and $\delta>0$.

Let us denote  by $\{ \varrho_r,\vu_r\}_{r\geq 1}$ a family of weak solutions to \eqref{usmiech_a3}. Then $\varrho_r \in L^{\beta+1} ((0,T)\times \Omega_r)$ and the continuity equation holds also in renormalized sense 
(provided $\varrho_r$, $\vu_r$ is extended by $0$ outside of $\Omega_r$). Finally  $(\varrho_r, \vu_r)$ satisfy for each $r \geq 1$ and $\delta>0$:
	\begin{equation}\label{p31_20}
	\| \varrho_r (t) \|_{L^1(\Omega_r)} = \| \varrho_r (0) \|_{L^1(\Omega_r)}\leq c \mbox{ for a.e. } t>0,
	\end{equation}
	\begin{equation}\label{p31_2}
	\sup_{t\in[0,T]} \| \varrho_r (t) \|^m_{L^m(\Omega_r)} \leq c(E_{\delta,r} (0),K),
	\end{equation}
	\begin{equation*}%\label{p31_31}
	\delta \sup_{t\in[0,T]} \| \varrho_r (t) \|^\beta_{L^\beta(\Omega_r)} 
	\leq c(E_{\delta,r} (0),K),
	\end{equation*}
	\begin{equation*}%\label{p31_3}
	\sup_{t\in[0,T]} \| \sqrt{\varrho_r (t)} \vu_r   \|^2_{L^2(\Omega_r)} 
	\leq c(E_{\delta,r} (0),K),
	\end{equation*}
	\begin{equation*}%\label{p31_4}
	\int_0^T    \|  \nabla_x \vu_r   \|^2_{L^2(\Omega_r)} \dt
	\leq c(E_{\delta,r} (0),K),
	\end{equation*}
	\begin{equation*}%\label{p31_7}
	\sup_{t\in[0,T]} \| \varrho_r  \Phi \|_{L^1(\Omega_r)}\leq c(E_{\delta,r} (0),K)  ,
	\end{equation*}
	\begin{equation}\label{p31_6}
	\sup_{t\in[0,T]} \| (K \ast \varrho_r) \varrho_r \|_{L^1(\Omega_r)}
	\leq c(E_{\delta,r} (0),K),
	\end{equation}
 and the
%\begin{equation*}
%\begin{split}
%	\frac{\rm d}{{\rm d}t} E(t)
%	+  \int_{\Omega_r}\left[\mu | \nabla_x \vu_r |^2 + (\lambda + \mu) | \Div \vu_r |^2 
%	+ \varrho_r |\vu_r|^2 \right] \dx \leq 0
%	\end{split}
%	\end{equation*}
%holds in $\D'(0,T)$, alongside with its 
integrated version of the energy inequality
\eqh{
E_{\delta,r}(t)+\int_0^t \int_{\Omega_r}\left[\mu | \nabla_x \vu_r |^2 + (\lambda + \mu) | \Div \vu_r |^2 
	+ \varrho_r |\vu_r|^2 \right]\dx\dt\leq E_{\delta,r}(0)
}
for $a.e.$ $t\in(0,T)$, where
\eqh{
	E_{\delta,r}(t)=\int_{\Omega_r} \left[ \frac{1}{2} \varrho_r | \vu_r |^2 + \frac{a}{m-1} \varrho_r^m  
	+ \frac{\delta}{\beta -1}\varrho^{\beta}_r + \right. & \left. \frac{1}{2}(K\ast\varrho_r) \varrho_r   + \varrho_r \Phi   \right] \dx,}
	\begin{equation*}%\label{energy0_delta_r}
	E_{\delta,r}(0)= 
	\int_{\Omega_r } \left[  
	\frac{1}{2} \frac{|\vm_{0,\delta}|^2}{\varrho_{0,\delta}}  + \frac{a}{m-1} \varrho_{0,\delta}^{m} 
	+ \frac{\delta}{\beta-1}\varrho_{0,\delta}^\beta + \frac{1}{2}(K \ast \varrho_{0,\delta})\varrho_{0,\delta} + \varrho_{0,\delta} \Phi 
	 \right] \dx.
	\end{equation*}
is bounded uniformly w.r.t. the parameters $r \geq 1$ and $\delta>0$. 
\end{Proposition}
 
 Let us notice that the boundeddnes of $E_{\delta,r} (0)$ results from the assumptions on the initial data and \eqref{ass_K}. Then the estimates \eqref{p31_2} --  \eqref{p31_6} are direct consequences of the energy inequality,  application of the Poincar\'e and the Young inequality, assumptions on $K$ (in particular see \eqref{p_ener}) and $\Phi$.  Moreover, \eqref{p31_20} is implied by the continuity equation. 
 
 \subsubsection{Existence of solutions on unbounded domain. Passing with $r  \to \infty$}\label{r_to_infty}

Let $\{ \Omega_r \}_{r>0}$ be given by \eqref{Omega_r}. 
Let us notice that we may extend $\varrho_r$  and $\vu_r$  by zero  on the whole $\R^3$ such that  Proposition~\ref{prop31} holds  true on $\R^3$.

Let $\varphi \in \D([0,T) \times \R^3)$ and $\vvarphi \in \D([0,T) \times \R^3; \R^3)$ be arbitrary, but fixed. Let us fix a number $\bar{r} \in (R, \infty)$ large enough, such that $\Omega_{\bar{r}}$ contains supports of both test functions. Thus, as in previous steps, we can pass to the limit in all terms of the weak formulation of the continuity and momentum equations \eqref{usmiech_a3} since our considerations can be reduced to the set $\Omega_{\bar{r}}$. 
 
The passage to the limit in the pressure term $p_r = a \varrho_r^m + \delta \varrho^\beta_r$, and in the nonlocal term $\vr_rK\ast\vr_r$, are similar as in the case of vanishing of the artificial pressure $\delta \to 0$, therefore the details of this discussion are postponed to the forthcoming subsection.

 \begin{Proposition}\label{prop_whole_space}
 Let \eqref{id_a1} be satisfied on whole $\R^3$ and let $0 < T < \infty$. Then there exists a bounded free energy weak solution $ (\varrho, \vu )$ to the problem
	 	\begin{align}
	\partial_t \varrho + \Div (\varrho \vu) 
	 = & \, 0 \quad \mbox{ in } (0,T) \times \R^3,\nonumber
	\\
	\partial_t(\varrho \vu) + \Div (\varrho \vu \otimes \vu ) + \nabla_x (a \varrho^m + \delta \varrho^\beta) 
	 = & \,\mu \Delta \vu + (\lambda + \mu) \nabla_x (\Div \vu) \label{usmiech_a4}
	 \\ &  - \varrho\vu - (\nabla_x K \ast \varrho)\varrho 
	- \nabla_x \Phi \varrho\quad \mbox{ in }  (0,T) \times \R^3 ,\nonumber
	\end{align}
for each fixed $\delta > 0$. In particular, the following energy inequality holds for a.e. $t\in(0,T)$
\eqh{
E_{\delta}(t)+\int_0^t \int_{\Omega_r}\left[\mu | \nabla_x \vu_r |^2 + (\lambda + \mu) | \Div \vu_r |^2 
	+ \varrho_r |\vu_r|^2 \right]\dx\dt\leq E_{\delta,r}(0),
}
where
$$
 	E_\delta (t) = \int_{\R^3} \left( \frac{1}{2} \varrho(t) |\vu (t) |^2 + \frac{a}{m-1} \varrho^m(t) + \frac{\delta}{\beta-1} \varrho^\beta + \frac{1}{2} (K \ast \varrho) \varrho + \Phi \varrho \right) \dx,
$$
	\begin{equation*}%\label{energy0_delta_r}
	E_{\delta}(0)= 
	\int_{\R^3 } \left[  
	\frac{1}{2} \frac{|\vm_{0,\delta}|^2}{\varrho_{0,\delta}}  + \frac{a}{m-1} \varrho_{0,\delta}^{m} 
	+ \frac{\delta}{\beta-1}\varrho_{0,\delta}^\beta + \frac{1}{2}(K \ast \varrho_{0,\delta})\varrho_{0,\delta} + \Phi \varrho_{0,\delta} 
	 \right] \dx.
	\end{equation*}
Denoting by $\{ \varrho_\delta, \vu_\delta \}_{\delta>0}$ solutions to \eqref{usmiech_a4}, the following estimates hold
	\begin{equation}\label{p431_2}
	\sup_{t\in[0,T]} \| \varrho_\delta (t) \|^m_{L^m(\R^3)} \leq c( E_{\delta}(0), K),
	\end{equation}
	\begin{equation*}%\label{p431_31}
	\delta \sup_{t\in[0,T]} \| \varrho_\delta (t) \|^\beta_{L^\beta(\R^3)} 
	\leq c( E_{\delta}(0), K),
	\end{equation*}
	\begin{equation}\label{p431_3}
	\sup_{t\in[0,T]} \| \sqrt{\varrho_\delta (t)} \vu_\delta   \|^2_{L^2(\R^3)} 
	\leq c( E_{\delta}(0), K),
	\end{equation}
	\begin{equation}\label{p431_4}
	\int_0^T    \|  \nabla_x \vu_\delta   \|^2_{L^2(\R^3)} \dt
	\leq c( E_{\delta}(0), K),
	\end{equation}
	\begin{equation*}%\label{p431_7}
	\sup_{t\in[0,T]} \| \varrho_\delta \Phi  \|_{L^1(\R^3)}\leq c E_{0,\delta} [\varrho_0,\vm_0, K]  ,
	\end{equation*}
and
	\begin{equation}\label{p431_6}
	\sup_{t\in[0,T]} \| (K \ast \varrho_\delta) \varrho_\delta \|_{L^1(\R^3)}\leq c( E_{\delta}(0), K).
	\end{equation}
%where 
%	\begin{equation*}
%	E_{0,\delta}[\varrho_0, \vm_0,K]= 
%	\int_{\R^3 } \left[  
%	\frac{1}{2} \frac{|\vm_{0,\delta}|^2}{\varrho_{0,\delta}}  + \frac{a}{m-1} \varrho_{0,\delta}^{m} 
%	+ \frac{\delta}{\beta-1}\varrho_{0,\delta}^\beta + \frac{1}{2}(K \ast \varrho_{0,\delta})\varrho_{0,\delta} + \Phi \varrho_{0,\delta}
%	 \right] \dx .
%	\end{equation*}
 \end{Proposition}
 
Let us notice again that estimates \eqref{p431_2} -- \eqref{p431_6} are direct consequences of the energy inequality.

Moreover the total mass is conserved in time.
\begin{Lemma}\label{lem_delta_lim}
Let $(\varrho_\delta, \vu_\delta)$ be a free energy weak solution to \eqref{usmiech_a4}.
Then the total mass is conserved in time
	\begin{equation}\label{p431_20}
	M_{\varrho_\delta}(t) = \int_{\R^3} \varrho_\delta(t, \cdot ) \dx  = \int_{\R^3} \varrho_{0,\delta} \dx.
	\end{equation}
\end{Lemma}
{\bf Proof.}
From energy estimates the confinement term provides 
	\begin{equation}\label{2nd_moment}
	\int_{\R^3} \Phi \varrho_{\delta,r} (t,x) \dx \leq C_T \quad \mbox{ for } t\in [0,T].
	\end{equation}
As \eqref{2nd_moment} and \eqref{ass_Phi_2} hold we infer by the  Chebyshev inequality   the following 
	\eq{\label{phi_33}
	\int_{\{ x \in \R^3: | x| \geq R\}} \varrho_{\delta,r}(t,x) \dx 
	&\leq
	\frac{1}{R^{1+\nu}} \int_{\R^3} |x|^{1+\nu}  \varrho_{\delta,r}(t,x) \dx\\
	 &\leq \frac{1}{R^{1+\nu}} \int_{\R^3} \Phi \varrho_{\varrho,r} (t,x) \dx \leq \frac{C_T}{R^{1+\nu}} \quad\mbox{ for } R>\bar{R}.
	}
Next, one deduce that
	\begin{equation}\label{mass_con_1}
	\int_{\R^3} \varrho_{\delta,r}(t) \dx = \int_{B(0,R)} \varrho_{\delta,r}(t) \dx + \int_{\R^ \setminus B(0,R)} \varrho_{\delta,r}(t) \dx 
	\leq  \int_{B(0,R)} \varrho_{\delta,r}(t) \dx  + \frac{C_T}{R^{1+\nu}}.
	\end{equation}
	Thus, by sending  $r \to \infty$ by weak convergence of density in $L^1_{loc}(\R^3)$  we deduce  for almost all $t \in (0,T)$ that 
 	$$M_{\delta,0} = \int_{\R^3} \varrho_{\delta,0} \dx =  \lim_{ r\to \infty} \int_{\R^3} \varrho_{\delta,r} (t) \dx =  \int_{B(0,R)} \varrho_\delta (t) \dx + O(\epsilon).$$          
Since the above holds for any $R$ we conclude that $M_{\varrho_\delta}(t)$ is a constant of motion as desired.	
\bigskip

Also the following higher local integrability result holds.

 	\begin{Lemma}\label{lem_theta}
	Let $\{ \varrho_\delta, \vu_\delta \}_{\delta>0}$ be a sequence of free energy weak solutions to the artificial pressure approximation \eqref{usmiech_a4} on $\R^3$. Then there exists a constant $c(T)$ independent of $\delta$ s.t.
	\begin{equation*}%\label{l_theta_1}
	\int_0^T \int_{B} a \varrho_\delta^{m + \theta} + \delta \varrho_\delta^{\beta+\theta} \dxdt \leq c(T,B) 
	\end{equation*}
for any compact $B\subset \R^3$, where $\theta = \min \{\tfrac{2}{3} m -1, \frac{1}{4}\}$.
	\end{Lemma}
 {\bf{Proof.}} 
 For the proof we follow \cite{CKT,FNP} with some slight modifications. For a given compact set $B$ we introduce a test function 
\eq{\label{lim_delta_3}
 	\vvarphi_\delta (t,x) =  \psi(t)  \A [\phi(x)  \varrho_\delta^{\theta}(t,x)]  \mbox{ where } \ \psi\in \D((0,T)),\  {\phi} \in \D(\R^3) ,\\
	 0\leq \phi \leq 1, \ \phi = 1 \mbox{ on }B, 
	 |\nabla_x \phi | \leq C \mbox{ on }\R^3 ,
	\mbox{ and } {\rm supp}\, \phi \subset 2B,
}
where $\A$ is given by \eqref{inv-div}. By \eqref{p31_2} and the choice of $\phi$ in \eqref{lim_delta_3}, we get
\eqh{
	\int_{\R^3} \phi(x) \varrho^\theta_\delta(t,x) \dx& \leq \int_{\{x: \varrho_\delta \leq 1\}} \phi(x) \varrho^\theta_\delta(t,x) \dx 
	+  \int_{\{x: \varrho_\delta > 1\}} \phi(x) \varrho^m_\delta(t,x) \dx \\
	&\leq |{\rm supp\,}\phi |  +  c(E_\delta(0),K) 
	 = c(\phi,\varrho_{0,\delta},\vm_{0,\delta})
}
and consequently by  \eqref{p431_2}, \eqref{p431_20} 
	\begin{equation*}%\label{theta_p_1}
	\| \phi \varrho^\theta_\delta \|_{L^\infty(0,T;L^1\cap L^{\frac{m}{\theta}}(\R^3))} \leq c(\phi,\varrho_{0,\delta}, \vm_{0,\delta}) .
	\end{equation*}
Then using Mikhlin multiplier theorem,  \eqref{Lizorkin}, the classical Sobolev embedding theorem we get (for more details see \cite[Section~ 3]{FP}):
	\begin{equation}\label{theta_p_2}
	\| \nabla_x   \vvarphi_\delta \|_{L^\infty(0,T; L^{r_1}(\R^3))} \leq c \quad \mbox{ with } r_1 \in \left(1, \frac{m}{\theta}\right],
	\end{equation}
	\begin{equation}\label{theta_p_3}
	\| \A [\Div (\phi \varrho_\delta^\theta \vu_\delta )] \|_{L^2(0,T; L^{r_2}(\R^3))} \leq c \quad \mbox{ with  }
	r_2 \in \left( 1, \frac{6m}{m + 6\theta} \right],
	\end{equation}
	\begin{equation}\label{theta_p_4}
	\| \vvarphi_\delta \|_{L^\infty( 0,T; L^{r_3}(\R^3))} \leq c \quad \mbox{ with } r_3 \in \left( \frac{3}{2}, \infty\right],
	\end{equation}
	\begin{equation}\label{theta_p_5}
	 \| \A [\varrho^\theta_\delta \vu_\delta \nabla_x \phi ] \|_{L^2( 0,T; L^{r_4}(\R^3))} \leq c  \quad \mbox{ with } r_4 \in \left( \frac{3}{2}, \infty \right),
	\end{equation}
and
	\begin{equation}\label{theta_p_6}
	 \| \A [\phi (\Div \vu_\delta) \varrho^\theta_\delta  ] \|_{L^2( 0,T; L^{r_5}(\R^3))} \leq c  
	 \quad \mbox{ with } r_5 \in\left( \frac{3}{2} , \frac{6m}{6\theta + m} \right].
	\end{equation}
The above estimates provide us that $\varphi_\delta$ is a proper test function for \eqref{usmiech_a4}$_2$ by density arguments,
and therefore we obtain
 	\begin{align}
	\int_0^T  \int_{\R^3} &  \left( a \varrho_\delta^{m+ \theta} + \delta \varrho_\delta^{\beta + \theta} \right) \phi \dxdt  \nonumber
	\\ =&\,   - \int_0^T \int_{\R^3}  (\varrho_\delta \vu_\delta) \partial_t \vvarphi_\delta \dxdt   
	- \int_0^T \int_{\R^3} (\varrho_\delta \vu_\delta \otimes \vu_\delta) : \nabla_x \vvarphi_\delta\dxdt
		\nonumber \\&
	+  \int_0^T \int_{\R^3}  \lr{ \mu \nabla_x \vu_\delta : \nabla_x \vvarphi_\delta 
	+ (\lambda + \mu) \Div \vu_\delta \Div \vvarphi_\delta}  \dxdt +  \int_0^T \int_{\R^3} (\nabla_x K \ast \varrho_\delta)\varrho_\delta \cdot \vvarphi_\delta \dxdt
	\nonumber \\&
	+  \int_0^T \int_{\R^3} \varrho_\delta \vu_\delta \cdot \vvarphi_\delta 
	+ \varrho_\delta \nabla_x \Phi \cdot  \vvarphi_\delta \dxdt 
	- \int_{\R^3} \vm_{0,\delta} \vvarphi_\delta(0,x) \dx 
	\nonumber \\ = &\, I_1 + I_2 + I_3 + I_4+ I_5.\label{lim_delta_5}
\end{align}
Then by \eqref{p431_20} and \eqref{p431_2}, the sequence $\{ \varrho_\delta\}_\delta$ is bounded in $L^\infty(0,T; L^{\tilde{m}}(\R^3))$ with $\tilde{m}=\frac{6m}{4m-6\theta}\leq m$.  This together with \eqref{p431_4} and the Sobolev embedding gives that $\|\varrho_\delta \vu_\delta\|_{L^2(0,T;L^p (\R^3))} \leq c$ with $p= \frac{6m}{5m - 6\theta}$.
Consequently, by   \eqref{theta_p_3}, \eqref{theta_p_5}, and \eqref{theta_p_6} for $p'= \frac{6m}{6\theta +m}$ we deduce
	\begin{align}
	|I_1|  =&  \left|  \int_0^T \int_{\R^3}  (\varrho_\delta \vu_\delta) \partial_t \vvarphi_\delta 
 \dxdt \right| \nonumber
	\\  
	\leq &
	\| \varrho_\delta \vu_\delta  \|_{L^2(0,T;L^{p}(\R^3))} c(T)
	\Big( \| \A[\phi (1-\theta) \varrho_\delta^\theta \Div \vu_\delta]\|_{L^2(0,T; L^{p'}(\R^3))} 
	+  \| \A [\varrho^\theta_\delta \vu_\delta \nabla_x \phi] \|_{{L^2(0,T; L^{p'}(\R^3))} }
	\nonumber \\ &\hspace{8cm} +    \| \A[ \Div (\phi \varrho^\theta_\delta \vu_\delta) \|_{L^2 (0,T; L^{p'}(\R^3))} \Big) \leq c.\label{lim_delta_6}
	\end{align}
Next by \eqref{theta_p_2}
	\begin{align}
 |I_2|  =  & \,\left|  \int_0^T \int_{\R^3}  
	(\varrho_\delta \vu_\delta \otimes \vu_\delta) : \nabla_x \vvarphi_\delta + \mu \nabla_x \vu_\delta : \nabla_x \vvarphi_\delta 
	+ (\lambda + \mu) \Div \vu_\delta \Div \vvarphi_\delta 
	 \dxdt \right| 
	\nonumber\\  
	\leq & \,
	\| \varrho_\delta \vu_\delta \|_{L^2(0,T;L^{p_2}(\R^3))} \| \vu\|_{L^2(0,T;L^{6}(\R^3))}
  \| \nabla_x  \vvarphi_\delta \|_{L^\infty(0,T; L^{p_3}(\R^3))} 
	\nonumber\\&
	+ c(T) \| \mu \nabla_x \vu_\delta \|_{L^2(0,T;L^{2}(\R^3))}
	  \| \nabla_x \vvarphi_\delta \|_{L^2(0,T; L^{2}(\R^3))} 
	\nonumber\\ &
	 + c(T) (\mu + \lambda) \| \nabla_x \vu_\delta \|_{L^2(0,T;L^{2}(\R^3))} 
	  \| \nabla_x   \vvarphi_\delta \|_{L^2(0,T; L^{2}(\R^3))}  \leq c ,\label{lim_delta_7}
\end{align}
where $p_2 = \frac{2m}{m+1}$ and $p_3$ is s.t. $\frac{1}{p_2} + \frac{1}{6} + \frac{1}{p_3} = 1$. 
For the estimate of the nonlocal term we use that due to Sobolev's imbedding we have
$$
\| \vvarphi_\delta \|_{L^\infty(0,T; L^{\infty}(\R^3))}\leq C\|\vr_\delta^\theta\|_{L^\infty(0,T; L^{m/\theta}(\R^3)}= C \|\vr_\delta\|_{L^\infty(0,T; L^m(\R^3))}^\theta,
$$
if only $m/\theta>3$, meaning that $m<3$. The case of $m\geq3$ is straightforward. We therefore can proceed analogously to \eqref{p_ener} to deduce
\begin{align*}%\label{p_bogo}
|I_3|=\left|\intTO{ (\nabla_x K \ast \varrho_\delta)\varrho_\delta \cdot \vvarphi_\delta } \right|
&\leq C \|\vr_\delta\|_{L^1(0,T; L^p(\R^3))}^2\|\Grad_x K\|_{L^{q}(\R^3)}\|\vr_\delta\|_{L^\infty(0,T; L^m(\R^3))}^\theta\nonumber\\
&\leq C \|\vr_\delta\|_{L^1(0,T; L^p(\R^3))}^2\|\Grad_x K\|_{L^q(\R^3)},
\end{align*}
where we denoted $q={\frac{p}{2(p-1)}}$. If $p=\frac{2q}{2q-1}\leq m$, i.e. $q\geq\frac{m}{ 2(m-1)}$, then we can interpolate $\|\vr_\delta (t)\|_{L^p(\R^3)}$ between $\|\vr_\delta(t)\|_{L^1(\R^3)}$ and $\|\vr_\delta(t)\|_{L^m(\R^3)}$ which are bounded uniformly in time. If on the other hand $p>m$, we can use the interpolation inequality between the spaces $L^\infty(0,T; L^m(\R^3))$ and $L^{m+\theta}((0,T)\times\R^3)$ to get
\eq{\label{p_bogo2}
 \|\vr_\delta(t)\|_{L^p(\R^3)}^2
\leq C\|\vr_\delta(t)\|_{L^m(\R^3)}^{2(1-\alpha)}\|\vr_\delta(t)\|_{L^{m+\theta}(\R^3)}^{2\alpha},
}
with the restrictions
\eqh{
m<\frac{2q}{2q-1},\quad \frac{1}{p}=\frac{1-\alpha}{m}+\frac{\alpha}{m+\theta},\quad 0\leq\alpha\leq1,\ \ \text{and}\quad 2\alpha<m+\theta.}
This leads to inequality $2\lr{1-\frac m p}<\theta$, or in other words $m<p<\frac{2m}{2-\theta}$, where $\theta=\min\{\frac23m-1,\frac14\}$, or finally, in terms of $q$ we get
$$
\frac{m}{2(m-1)+\theta}<q,  \quad \theta=\min\left\{\frac23m-1,\frac14\right\}.
$$
Then to get rid of the term $\|\vr_\delta(t)\|_{L^{m+\theta}(\R^3)}^{2\alpha}$ from \eqref{p_bogo2}  that contributes to the RHS of \eqref{lim_delta_5} 
we use the Young inequality.
Next, by \eqref{theta_p_4}
 	\begin{align}
	| I_4 & + I_5 | =   \left|  \int_0^T \int_{\R^3} \varrho_\delta \vu_\delta \cdot \vvarphi_\delta 
	+ \nabla_x \Phi \varrho_\delta \cdot \vvarphi_\delta  \dxdt
	- \int_{\R^3} \vm_{0,\delta} \vvarphi_\delta(0,x) \dx   \right| 
	\nonumber\\
	 \leq & 
	 \| \varrho_\delta \vu_\delta  \|_{L^1(0,T; L^p(\R^3))}  \| \vvarphi_\delta \|_{L^\infty(0,T; L^{p'}(\R^3))} 
	\nonumber \\ & 
	+
	 \left(\| \nabla_x \Phi \|_{L^{m'}((0,T)\times (\Omega_r \cap B_{R_o}))} \| \varrho_\ep\|_{L^m((0,T)\times (\Omega_r \cap B_{R_o}))} + \| \Phi \varrho_\ep\|_{L^1((0,T) \times (\Omega_r  \setminus B_{R_o}))}\right) \| \vvarphi_\delta \|_{L^\infty(0,T; L^{\infty}(\R^3))}
	\nonumber \\&
	+ \| \vm_{0,\delta}\|_{L^1}   \|  \A [\phi \varrho_{0,\delta}^{\theta}] \|_{L^\infty(\R^3)}  
	\leq c . \label{lim_delta_9}
	\end{align}
Here, $p= \frac{6m}{5m - 6\theta}$ and $B_{R_o}$ is a ball centred at zero with radius $R_o$  given by \eqref{ass_Phi_1}.

Collecting the estimates in \eqref{lim_delta_5}, \eqref{lim_delta_6}, \eqref{lim_delta_7}, and \eqref{lim_delta_9}, we infer the claim in Lemma~\ref{lem_theta}. 
 
 \subsubsection{Artificial pressure limit $\delta \to 0$}
 In previous subsections we provided existence of weak solutions to the artificial pressure approximation on $\R^3$. Let us now pass with $\delta \to 0$. For this reason let us consider general initial data $\varrho_0,\ \vm_0$ satisfying \eqref{initial_data} and \eqref{com_con}. Then, we can find a sequence of approximations to the initial data as
	\begin{equation*}%\label{delta_limit_id}
	\varrho_{0,\delta} \in C^{2+ \nu}(\R^3) \mbox{ s.t. }
	0<\varrho_{0,\delta}(x) \leq \delta^{-\frac{1}{\beta}}, \quad \varrho_{0,\delta} \to \varrho_0 \mbox{ in }L^m(\R^3) \cap L^1_{+}(\R^3)
	\mbox{ as }\delta \to 0
	\end{equation*}
and 
$$
	\vm_{0,\delta} \to \vm_0 \quad\mbox{ in }L^1(\R^3)\mbox{ as }\delta \to 0.
$$ 
With \eqref{usmiech_a4}$_1$, \eqref{p431_20},  \eqref{p431_2}, and \eqref{p431_3} at hand, we can deduce, similarly  as for \eqref{rho_ep_weak_CLm}, that
 	\begin{equation}\label{rho_delta_weak_CLm}
	\varrho_\delta \to \varrho \quad \mbox{ in }C_{weak}([0,T];L^m(\R^3)).
	\end{equation}
Similarly, by \eqref{usmiech_a4}$_2$ and the estimates obtained in Proposition~\ref{prop_whole_space}, we infer
 	\begin{equation*}%\label{rhou_delta_weak_CLm}
	\varrho_\delta \vu_\delta \to \varrho \vu \quad \mbox{ in }C_{weak}([0,T];L^{\frac{2m}{m+1}}(\R^3)).
	\end{equation*}
Then, using \eqref{p431_4} and Div-Curl type argument we deduce also that
 	\begin{equation*}%\label{rhou_delta_weak_CLm}
	\varrho_\delta \vu_\delta\otimes\vu_\delta \weak \varrho \vu \otimes\vu\quad \mbox{ in }L^2(0,T;L^{\frac{6m}{4m+3}}_{loc}(\R^3)).
	\end{equation*}
Due to Lemma~\ref{lem_theta} we also have
  $$
  \varrho^m_\delta \rightharpoonup \Ov{ p}= \overline{\varrho^m} \quad \mbox{ in } L^{\frac{m +  \theta }{m}}_{loc}((0,T)\times \R^3).
  $$ 

Therefore, the estimates given in Proposition~\ref{prop_whole_space}, Lemma~\ref{lem_theta} and the Banach-Alaoglu theorem allow us to pass to the limit as $\delta \to 0$ in the continuity and momentum equation \eqref{usmiech_a4}. We find that
 	\begin{equation*}%\label{lim_delta_1}
	\int_{0}^T \int_{\R^3} \varrho \partial_t \varphi - (\varrho \vu ) \cdot \nabla_x \varphi \dxdt = 
	\int_{\R^3} \varrho_0 \varphi(0,x) \dx
	\end{equation*}
 for all $\varphi \in \D( [0,T) \times \R^3)$ where the continuity equation is also satisfied in renormalized sense, i.e.  \eqref{weak_1} is satisfied, and 
 	\begin{equation*}%\label{lim_delta_2} 
	\begin{split}
	\int_0^T \int_{\R^3} & (\varrho \vu) \partial_t\vvarphi + (\varrho \vu \otimes \vu) : \nabla_x \vvarphi 
	+ \overline{p} \Div \vvarphi \dxdt
	\\  = & \int_{0}^{T} \int_{\R^3}  \mu \nabla_x \vu : \nabla_x \vvarphi + (\lambda + \nu) \Div \vu \Div \vvarphi 
	+ (\varrho \vu) \vvarphi  + \left(\overline{(\nabla_x K \ast \varrho)\varrho} +\varrho   \nabla_x\Phi  \right) \cdot \vvarphi \dxdt 
	\\ & - \int_{\R^3} \vm_0 \cdot \vvarphi(0,x) \dx 
	\end{split}
	\end{equation*}
 for all $\vvarphi \in \D( [0,T) \times \R^3)$. 
 
The terms with the bars  cannot be identified yet, as we do not know whether the density sequence $\vr_\delta$ converges strongly to $\vr$. In order to prove it, we use Feireisl's technique \cite{F_com} that allows to treat the cases of non square-integrable densities. We introduce the following truncation operator:
\eqh{
T_k(z) = k T\Big(\frac{z}{k}\Big), \quad z \in \R,\quad k \ge 1,
}
with $T \in C^\infty (\R) $ such that
\begin{equation*}
T(z) = z~\text{for}~z\le 1, \quad T(z) =2~\text{for}~z \ge 3,\quad T~\text{concave, non-decreasing}.
\end{equation*}
We use as the test function  in the approximate momentum equation
the function
\[
 \vvarphi_\delta= \psi(t) \nabla(x) \Delta^{-1}[\phi T_k(\vr_\delta)] = \psi(t) \phi(x) \A [ T_k(\vr_\delta)],\,k\in \mathbb{N},
\]
and for the limit equation
the test function
\[
 \vvarphi=\psi(t) \phi(x)\nabla \Delta^{-1}[\Ov{T_k(\vr)}] = \psi(t) \phi(x) \A [\Ov{T_k(\vr)}] ,\,k\in \mathbb{N}.
\]
Here, $\psi \in \D((0,T))$ and $\phi \in \D(\R^3)$. 
Subtracting the resulting expressions, we obtain the following
\eq{ \label{ws28a}
&\lim_{\delta \to 0^+} \int_0^T \psi \int_\Omega \phi \Big( a\vr_\delta^m T_k(\vr_\delta) - (\lambda + 2 \mu) \Div \vu_\delta T_k(\vr_\delta) \Big)
\, \dxdt
\\
&\quad-\int_0^T \psi \int_\Omega \phi \Big( a\Ov{\vr^m} \ \Ov{T_k(\vr)} - (\lambda + 2 \mu) \Div\vu \Ov{T_k(\vr)} \Big)
\, \dxdt
\\
&=-\lim_{\delta \to 0^+} \int_0^T \psi \int_\Omega \phi \vr_\delta\Grad K\ast\vr_\delta\cdot  \A[ T_k(\vr_\delta)]\, \dxdt
+ \int_0^T \psi \int_\Omega \phi \Ov{\vr\Grad K\ast\vr}\cdot  \A [ \Ov{T_k(\vr)}]\, \dxdt\\
&\quad+ \lim_{\delta \to 0^+}{\mathcal R}(\vr_\delta,\vu_\delta)-{\mathcal R}(\vr,\vu).}
The last two terms contain all the contributions coming from the time derivative of momentum, the acceleration, confinement, and friction. Following \cite{F_com} and \cite{CKT}, we can show that these two terms cancel out when $\delta\to 0$.
The other two terms on the r.h.s. of \eqref{ws28a} cancel out as well. Indeed, 
according to our definition, $T_k(\vr_\delta)$ is a good renormalization function, we therefore have
\begin{equation*} %\label{8.2bbb}
\partial_t(T_k(\vr_\delta))+\Div(T_k(\vr_\delta)\vu_\delta)+(\vr_\delta T'_k(\vr_\delta)-T_k(\vr_\delta))\Div\vu_\delta=0,
\end{equation*}
in the sense of distributions. This provides and estimate of the time derivative of $T_k(\vr_\delta)$, and since the operator $\A$ ''gains'' one spatial derivative, we find that
\[
 \A [ T_k(\vr_\delta) ] \to  \A[ \overline{T_k(\vr)}] \ \mbox{in}\ C([0,T] \times\R^3).
\]
On the other hand, we know that
\[ \vr_\delta\Grad K\ast\vr_\delta\to\Ov{\vr\Grad K\ast\vr}\ \ \mbox{weakly in}\ L^p((0,T) \times \R^3),\quad  \mbox{for some}\  p>1,\]
and therefore 
\eqh{
\lim_{\delta \to 0^+} \int_0^T \psi \int_\Omega \phi \vr_\delta\Grad K\ast\vr_\delta\cdot  \A[ T_k(\vr_\delta)]\, \dxdt= \int_0^T \psi \int_\Omega \phi \Ov{\vr\Grad K\ast\vr}\cdot  \A [ \Ov{T_k(\vr)}]\, \dxdt.
}
As a consequence we obtain the so-called effective viscous flux equality
\eqh{
&\lim_{\delta \to 0^+} \int_0^T \psi \int_\Omega \phi \Big( \vr_\delta^m T_k(\vr_\delta) - (\lambda + 2 \mu) \Div \vu_\delta T_k(\vr_\delta) \Big)
\, \dxdt
\\
&\quad=\int_0^T \psi \int_\Omega \phi \Big( \Ov{\vr^m} \ \Ov{T_k(\vr)} - (\lambda + 2 \mu) \Div\vu \Ov{T_k(\vr)} \Big)
\, \dxdt.
}
Once this equality is guaranteed, the proof of strong convergence of the density follows exactly as in  \cite{FNP} . The  slight modification concerning additional damping, confinement, and nonlocal terms similar to the case of $\ep \to 0$ limit above.
Notice that when passing to the limit with $\delta \to 0$ in our weak formulation due to the compact support of the test functions, we may restrict our considerations to some $\Omega_{\bar{r}}$ set containing the support of a test function. In particular 
\eqh{
 \delta \varrho^{\beta}_\delta \to 0 \mbox{ in }L^1(0,T; L^1_{loc}(\R^3), \qquad 
\qquad \varrho_\delta \to \varrho \quad \mbox{ in } L^m(0,T;L^m_{loc}(\R^3)),}
and so
\eqh{
\Ov{p(\vr)}=a\vr^m,\quad  \overline{(\nabla_x K \ast \varrho)\varrho}=(\nabla_x K \ast \varrho)\varrho,\quad \mbox{ a.e. in }\R^3.
}
With this Theorem~\ref{t_existence} is proved. 

\subsection{Properties of solutions}
\begin{Corollary}
Let $(\varrho, \vu)$ be a free energy weak solution to \eqref{main_system}.
Then the total mass is conserved in time
	\begin{equation}\label{mass_cons}
	M_\varrho(t) = \int_{\R^3} \varrho(t, \cdot ) \dx  = const .
	\end{equation}
\end{Corollary}
The above is a consequence of Lemma~\ref{lem_delta_lim}.
\begin{Lemma}\label{first_moment_bdd}
Let assumptions of Theorem~{\rm\ref{t_existence}} be satisfied. Let us assume that  $x \varrho_0$ is bounded in $L^1(\R^3)$. Then $x \varrho (t)$ is bounded in $L^1(\R^3)$ for all times and it solves the equation
$$
\frac{{\rm d^2}}{{\rm d}t^2} \int_{\R^3} x \varrho \dx = \frac{{\rm d}}{{\rm d}t} \int_{\R^3} \varrho \vu \dx = - \int_{\R^3}  \left( \varrho \vu +\varrho \nabla_x \Phi  \right) \dx\,.
$$
\end{Lemma}	
{\bf Proof.} Let us multiply the continuity equation of the system \eqref{usmiech_a3}$_1$ by 
$\psi_\zeta \in C^\infty_c (\R^3)$ such that $\psi_\zeta (x) \to  x $ in $L^{m'} \cap  L^{\frac{2m}{m-1}}(\Omega_r)$ as $\zeta \to 0$ and let us integrate over $\R^3$ (extending $\varrho_{\delta,r}$ and $\vu_{\delta,r}$ by zero on $\R^3 \setminus \Omega_r$). Letting $\zeta \to 0$, $r\to \infty$ and $\delta \to 0$ we obtain 
	\begin{equation}\label{conf_1}
	\frac{{\rm d}}{{\rm d}t} \int_{\R^3} x \varrho \dx = \int_{\R^3} \varrho \vu \dx\,,
	\end{equation}
as $\varrho(t) \in L^m(\R^3)$ and $\varrho \vu (t) \in L^{\frac{2m}{m+1}}(\R^3)$ for a.e. $t\in (0,T)$. 
When passing to the limit the first moment is controlled in a similar way to \eqref{phi_33} due to assumption \eqref{ass_Phi_2}. Now, we repeat an analogous procedure to momentum equation \eqref{usmiech_a3}$_2$ with  $\psi = 1$. Note that the boundary term coming from the pressure function and nonlocal term vanish due to the symmetry of $K$ and the proper choice of boundary conditions for the approximation. Consequently, we get
	\begin{equation*}\label{conf_2}
	\begin{split}
	\frac{{\rm d}}{{\rm d}t} \int_{\R^3} \varrho \vu \dx & = - \int_{\R^3} \left( (\nabla_x K \ast \varrho) \varrho + \varrho \vu + \varrho \nabla_x \Phi  \right) \dx
	\\ & =  - \int_{\R^3}  \left( \varrho \vu +\varrho  \nabla_x \Phi  \right) \dx\,,
	\end{split}
	\end{equation*}
and solving the ODE, we conclude that
	\begin{equation}\label{conf_3}
	\int_{\R^3} (\varrho \vu)(t) \dx = e^{-t}\left( \int_{\R^3} \vm_0\dx   
	- \int _0^t  e^s \int_{\R^3}  \varrho\nabla_x \Phi \dx\, {\rm d}s \right).
	\end{equation}
Therefore \eqref{conf_1} together with \eqref{conf_3} give 
  \begin{equation}\label{1mom}
  \int_{\R^3} (x \varrho)(\tau) \dx  = \int_0^\tau e^{-t} \dt  \left(  \int_{\R^3} \vm_0\dx - \int _0^t   e^s \int_{\R^3}   \varrho \nabla_x \Phi \dx\, {\rm d}s   \right)+\int_{\R^3} x \varrho_0 \dx.
  \end{equation}
  In order to show that the RHS of \eqref{1mom} is bounded we notice that
  $$
  \left| \int_{\R^3}  \varrho \nabla_x \Phi  \dx \right| \leq  \int_{\R^3} | \varrho  \nabla_x \Phi | \dx
  \leq \int_{|x|\leq R}  \varrho \nabla_x \Phi  \dx + \int_{|x| > R}  \varrho  \Phi  \dx
  $$
  and the RHS of the above is bounded as $\Phi \in W^{1,\infty}_{loc}(\R^3)$, \eqref{mass_cons}, and $\Phi \varrho \in L^\infty(0,T; L^1(\R^3))$ by energy estimates.
This finishes the proof of Lemma~\ref{first_moment_bdd}. 
  
\subsubsection{Global in time existence}
 
 We prove global in time existence by patching local-in-time solutions since the estimates on the energy provide uniform in time bounds. 
 
 \begin{Lemma}\label{L:global}
 Let $\Omega=\R^3$ or let $\Omega$ be smooth bounded set. Let assumptions of 
 Theorem~{\rm\ref{t_existence}} or Theorem~{\rm\ref{t_existenceb}} be satisfied respectively.
 Then solutions to the system \eqref{main_system} are global in time, namely exist on time interval $[0,\infty)$.
 \end{Lemma}
 
Fix any $T_0 >0$. Then there exists weak solution $(\varrho_1, \vu_1)$ on $[0, 2T_0)$. Let us introduce 
	\begin{equation*}
	\zeta^\kappa(t) = \left\{
\begin{array}{ccc}
 1 &    t\in[0,T_0- \kappa]  \\
 -\frac{1}{\kappa} |T_0 - t| & t \in (T_0 - \kappa, T_0)     \\
  0 &    \mbox{ otherwise}.
\end{array}
\right.
	\end{equation*}
Let us set $\varphi \zeta^\kappa$ with $\varphi \in \D( [0,T) \times \Omega)$ as a test function for continuity equation (in the renormalized sense) 	
	\begin{equation*}
	\begin{split}
	-\int_{T_0 - \kappa}^{T_0}  \int_{\Omega} b(\varrho_1) \varphi \kappa^{-1} \dxdt 
	& + \int_{0}^{T_0} \int_\Omega\left(b(\varrho_1) \partial_t \varphi \zeta^\kappa +  b(\varrho_1) \vu_1 \cdot \nabla \varphi \zeta^\kappa - (b'(\varrho_1) \varrho_1 - b(\varrho_1)) \Div \vu_1 \varphi \zeta^\kappa \right)  \dxdt   
	\\ &
	= -\int_{\Omega} b(\varrho_0)  \varphi (0, \cdot) \zeta^\kappa(0) \dx
	\end{split}
	\end{equation*}
for all $\varphi \in \D((0,T)\times \R^3)$. Letting $\kappa \to 0$ 
	\begin{equation*}
	\begin{split}
	\int_0^{T_0} \int_{\Omega}& \left( b(\varrho_1) \partial_t \varphi 
	+ (b(\varrho_1)  \vu_1) \cdot \nabla_x \varphi 
	- (b'(\varrho_1)\varrho_1  - b(\varrho_1) ) \Div \vu_1 \varphi  \right)
	\dxdt 
	\\ &
	= 
	 \int_{\Omega} b(\varrho (T_0, \cdot)) \varphi (T_0, \cdot) \dx - \int_{\Omega} b(\varrho_0) \varphi(0,\cdot) \dx
	 \end{split}
	\end{equation*}

Hence $(\varrho_1, \vu_1)$ is a free energy weak solution of the continuity equations on the closed interval $[0, T_0]$ with extra boundary term at $t=T_0$. We may apply the same argument to the momentum equation so that $(\varrho_1, \vu_1)$ is a weak solution on $[0, T_0]$ with boundary term at time $t=T_0$.
By uniform in time bounds at $T_0$ one can construct new solutions $(\varrho_2, \vu_2)$ defined on $[T_0,2T_0) \times \Omega$ such that
	$$E[\varrho_2,\vu_2] \leq E[\varrho_0,\vu_0] \quad t \in [T_0, 2T_0].$$
Then the couple $(\varrho, \vu)$ given by 
	\begin{equation*}
	(\varrho,\vu)(t) = 
	 \left\{
	\begin{array}{ccc}
	(\varrho_1,\vu_1) & t\in[0,T_0] \\
	(\varrho_2,\vu_2) & t\in(T_0,2T_0]
	\end{array}
	\right.
	\end{equation*}
is a solution on $[T_0,2T_0) \times \Omega$. A solution for all times is readily obtained by iterating the above procedure.
				
%%%%%%%%%%%%%%%%%%%%%%%%%%
	
\section{Long time asymptotics}\label{s:long_time}
In this section we consider the long-time asymptotic of solutions to \eqref{main_system}. In particular, we give here a proof of Theorem~\ref{t_longtime}. Let $\{ t_n \}_{n\geq 1}$ be a sequence s.t. $t_n \to \infty $ as $n\to \infty$ and let us define the sequences: 
	\begin{equation}\label{d:n-seq}
	\varrho_n (t,x) = \varrho (t+ t_n, x), \quad\quad \vu_n (t,x) = \vu (t+ t_n, x)  \quad \mbox{ for } t\in(-1,2), \ x\in \Omega \,.
	\end{equation}
Let us note that for each $n \in \mathbb{N}$ a couple $(\varrho_n, \vu_n)$ is a weak solution to the system \eqref{main_system} in a sense of Definition~\ref{defsol} and obtained in Theorem~\ref{t_existence}, Theorem~\ref{t_existenceb} respectively on $\Omega=\R^3$ or $\Omega$ bounded. Using the bounds obtained in previous sections we get the following estimates.
 
	\begin{Lemma}\label{long_time}
	Let $\Omega = \R^3$ or be a bounded smooth domain in $\R^3$. Let $(\varrho_n, \vu_n)$ be given by \eqref{d:n-seq} where $(\varrho, \vu)$ is a solution to \eqref{main_system} given by Theorem~{\rm\ref{t_existence}}, Theorem~{\rm\ref{t_existenceb}} respectively. Then
	\begin{equation}\label{Lem_LT_1}
	\lim_{n \to \infty} \int_{- 1}^{ 2} \| \nabla_x \vu_n \|_{L^{2}(\Omega)} 
	+ \| \varrho_n |\vu_n|^2 \|_{L^{1}(\Omega)}  \dt = 0,
	\end{equation}
	\begin{equation}\label{Lem_LT_3}
	\sup_{n\in N,\ t\in(-1,2) }  \left( \| \varrho_n\|_{L^m(\Omega))} +  \| \varrho_n (\nabla_x K \ast \varrho_n)\|_{L^1(\Omega)} + \| \varrho_n \Phi  \|_{L^1(\Omega)} 
	+ \| \varrho_n |\vu_n|^2 \|_{L^1(\Omega)} \right) < c,
	\end{equation}
	\begin{equation}\label{Lem_LT_4}
	\int_{\Omega} \varrho_n (t) \dx = \int_{\Omega} \varrho_0 \dx, 
	\end{equation}
	and 
	\begin{equation}\label{Lem_LT_2}
	 \int_{1}^{ 2} \int_{\R^3} \varrho_n^{m+ \theta} \phi \dxdt \leq c (B) \quad \mbox{ for all } \tau >0 
	 \mbox{ and }\phi \in C^\infty_c(\R^3), \ 0\leq \phi \leq 1, \ {\rm supp\, }\phi \subset B
	\end{equation}
where  $\theta = \min \{ \frac{2}{3} m - 1, \frac{1}{4} \}$, here $\varrho_n$ is extended to be zero on $\R^3\setminus \Omega$ if needed.
	\end{Lemma}
First, let us notice that \eqref{u_to_zero_a} and \eqref{u_to_zero_b} are direct consequences of \eqref{Lem_LT_1}.
Moreover due to \eqref{Lem_LT_1}$_2$,  \eqref{Lem_LT_3}$_1$
	\begin{equation}\label{vvv}
	\| \varrho_n \vu_n \|_{L^2(-1,2;L^{\frac{2m}{m+1}}(\Omega))} \leq \| \sqrt{\varrho_n} \vu_n \|_{L^2(-1,2;L^2(\Omega))} \| \sqrt{\varrho_n}\|_{L^\infty(-1,2;L^m(\Omega))} \to 0 \quad\mbox{ as } n\to \infty ,
	\end{equation}
what comes from the H\"older inequality for Bochner spaces.

The compactness deduced from the uniform in time energy bounds provide the existence of functions $\varrho_s$, $\overline{p}$ s.t. 
	\begin{equation}\label{n_to_s_1}
	\varrho_n \weak \varrho_s \quad \mbox{ weakly in }L^m( (-1,2) \times \Omega),
	\end{equation}
	\begin{equation}\label{n_to_s_2}
	\begin{split}
	a\varrho_n^m = p_n \weak \overline{p} \mbox{ weakly in } & L^{q}( -1,2; L^q_{loc} (\R^3)) \mbox{ or weakly in }L^{q}( -1,2; L^q (\Omega)) \mbox{ for }\Omega \mbox{ bounded }
	\\ &  \mbox{ with }1< q \leq \frac{m+ \theta}{m}.
	\end{split}
	\end{equation}
Taking the limit $n \to \infty $ in continuity equation, by \eqref{vvv} and \eqref{n_to_s_1} we get
	\begin{equation*}
	\int_{-1}^{2} \int_{\R^3} \varrho_s \varphi_t \dxdt = 0 \quad \mbox{ for all } \varphi_t \in C_c^\infty ((-1,2) \times \Omega) .
	\end{equation*}
Thus $\varrho_s$ is independent of time.	
Notice that similarly as for \eqref{mass_con_1} for $\Omega=\R^3$ we have
	\begin{equation*}
	\begin{split}
	\int_{\R^3} \varrho_0 \dx & = \lim_{n \to \infty } \frac{1}{3} \int_{-1}^{2} \int_{\R^3} \varrho_n \dxdt =
	\lim_{n\to\infty} \left[ 
	\frac{1}{3} \int_{-1}^{2} \int_{B(0,R)} \varrho_n \dxdt  + \frac{1}{3} \int_{-1}^{2} \int_{\R^3\setminus B(0,R)} \varrho_n \dxdt 
	\right]
	\\ &
	=  \int_{B(0,R)} \varrho_s \dxdt  + \frac{1}{3} 3 \frac{C_T}{R^{1+\nu}}.
	\end{split}
	\end{equation*}
Then for any $\ep$ we find $R$ s.t $ \frac{C_T}{R^2} < \ep$. Therefore
$\int_{\R^3} \varrho_s \dx = \int_{\R^3} \varrho_0 \dx  .$
Observe that for the case $\Omega$ bounded we do not need to use assumption \eqref{ass_Phi_2} in order to prove an analogous property. Therefore we can conclude that 
	$$\int_{\Omega} \varrho_s \dx = \int_{\Omega} \varrho_0 \dx  .$$

Next we pass to the limit $n\to \infty$ in the momentum equation.
Directly by \eqref{Lem_LT_1}, \eqref{vvv} we have that
	\begin{equation}\label{n_to_s_3}
	\begin{split}
	\lim_{n\to \infty}
	\int_{-1}^{2}\int_{\Omega}&  \left( \varrho_n \vu_n \cdot \partial_t\vvarphi  + (\varrho_n \vu_n \otimes \vu_n  - \mu \nabla_x \vu) : \nabla_x \vvarphi   - (\lambda + \mu) \Div \vu \Div \vvarphi \right) \dxdt 
	\\ & +\int_{\Omega} (\varrho_n\vu_n)(-1) \cdot \varphi(-1) \dx=0
	\end{split}
	\end{equation}
for any $\vvarphi \in \D ([-1,2) \times \Omega ;\R^3)$. Then proceeding as in as in Section~\ref{r_to_infty} %(see \eqref{K_lim_delta}, similarly to \eqref{K_lim_ep}) 
and by 
\eqref{n_to_s_1}, $\Phi \in W^{1,\infty}_{loc}(\R^3))$, \eqref{n_to_s_2} we infer
	\begin{equation}\label{n_to_s_4}
	\begin{split}
	\lim_{n \to \infty }&  \int_{-1}^{2}\int_{\Omega} \left( a\varrho_n^m \Div   \vvarphi  - (\nabla_x K \ast \varrho_n )\varrho_n \cdot \vvarphi - \varrho_n  \nabla_x \Phi \cdot \vvarphi \right) \dxdt 
	\\ & = \int_{-1}^{2}\int_{\Omega} \left(  \overline{p} \Div   \vvarphi  - \Ov{(\nabla_x K \ast \varrho_s )\varrho_s} \cdot \vvarphi - \varrho_s \nabla_x \Phi  \cdot \vvarphi \right) \dxdt 
	\end{split}
	\end{equation}
for any $\vvarphi \in \D ([-1,2) \times \Omega ;\R^3)$. Summarising \eqref{n_to_s_3} and \eqref{n_to_s_4} since   $\varrho_s$ is independent of time we obtain the balance of forces relation
	\begin{equation*}
	\nabla_x \overline{p} = - \Ov{(\nabla_x K \ast \varrho_s)\varrho_s }  - \nabla_x \Phi \varrho_s \quad\mbox{ in }\D'(\Omega).
	\end{equation*}
It only left to prove that 
\begin{equation}\label{p_K_n_lim}
\overline{p} =  a \varrho^m_s \quad \mbox{ and } \quad \Ov{(\nabla_x K \ast \varrho_s)\varrho_s}  = (\nabla_x K \ast \varrho_s)\varrho_s.
\end{equation}
In particular we claim that
	\begin{equation*}%\label{rho_n_s_q_loc}
	 \varrho_n \to \varrho_s \quad \mbox{ strongly in }L^q((-1,2)\times B) \quad\mbox{ for } q\in[1,m)
	 \mbox{ for any compact }B\subset \R^3 .
	 \end{equation*}
Since \eqref{Lem_LT_3} holds, it is enough to show that 
	\begin{equation}\label{rho_n_in_L1}
	\varrho_n \to \varrho_s \quad \mbox{ in } L^1((-1,2) \times B ) .
	\end{equation}
With this aim, we essentially  follow the proof of \cite[Theorem~1.1]{F_com}. 
Uniform in $n$ bounds in Lemma~\ref{long_time} imply that we can adopt the arguments from \cite{F_com} to show that 
	\begin{equation}\label{rho_n_s_ae}
	\varrho_n \to \varrho_s \quad \mbox{ a.e. in }(-1,2) \times\Omega
	\end{equation}
	and consequently \eqref{rho_n_in_L1} with \eqref{p_K_n_lim}.
Moreover let us notice that \eqref{rho_n_s_ae}, \eqref{Lem_LT_3}, \eqref{Lem_LT_4} and \eqref{Lem_LT_2} give uniform integrability of $\{ \varrho_n\}_{n\geq 1}$ in $L^q(\Omega)$ with $q\in [1,m]$ (in both cases $\Omega$ bounded and unbounded) and consequently by Vitali's Theorem 
 	\begin{equation*}%\label{T_t_m}
	\varrho_n \to \varrho_s \quad \mbox{ strongly in } L^m((-1,2)\times\Omega)
	\end{equation*}
	 holds true.
Moreover as for \eqref{rho_delta_weak_CLm}, we get 
	 \begin{equation*}%\label{time_ae_Lm}
	\varrho_n \to \varrho_s \quad \mbox{ in } C_{weak} ([-1,2];L^m(\Omega)).
	\end{equation*}
Consequently we infer that
	$$\varrho(t_n) \to \varrho_s \quad \mbox{ strongly in }L^m(\Omega) \mbox{ as } t_n\to \infty.$$
This finishes the proof of the characterization of the $\omega$-limit set in $L^m(\Omega)$ and of  in Theorem~\ref{t_longtime}, since if we are in the whole space case the confinement term allows us to pass to the limit the first moment if it has a limiting value.

Finally, Lemma~\ref{first_moment_bdd} gives us two particular cases in which we can discuss the long time asymptotics of the first moment. If we have quadratic confinement  Lemma~\ref{first_moment_bdd} gives a closed second order differential equation for the first moment. Moreover, the confinement provided by the energy estimate on the integral of $\Phi\varrho$ allows to pass to the limit the first moment.
  
\begin{Corollary}\label{first_moment_whole2}
Let assumptions of Theorem~{\rm\ref{t_existence}} be satisfied and let $\Phi =\tfrac{|x|^2}{2}$. Then we have that 
$$
\lim_{t\to\infty}  \int_{\R^3} x \varrho(t,x) \dx  = 0.
$$
Moreover, the $\omega$-limit set $\omega(\varrho)$ of the solution in $L^m(\Omega)$ consists of densities centered at zero, $c_0=0$ in Theorem~{\rm\ref{t_longtime}}.
\end{Corollary}  

If we are in the bounded and symmetric domain case, we can allow to have no confinement and also get a nice behavior asymptotically of the first moment.
  
\begin{Corollary}\label{first_moment_bdd_bdd}
Let assumptions of Theorem~{\rm\ref{t_existenceb}} be satisfied and let $\Phi \equiv 0$ and $\Omega$ be symmetric, i.e. $\Omega=-\Omega$. Then we have that 
$$
  \int_\Omega x \varrho(t,x) \dx  = \int_\Omega x \varrho_0(x)\dx + (1-e^{-t}) \int_\Omega \vm_0 \dx.
$$
Moreover, the $\omega$-limit set $\omega(\varrho)$ of the solution in $L^m(\Omega)$ consists of densities centered at the asymptotic value
$$
c_0=\int_\Omega x \varrho_0(x)\dx + \int_\Omega \vm_0 \dx\,
$$
in Theorem~{\rm\ref{t_longtime}}.
\end{Corollary}
  
Notice that the symmetry of the domain $\Omega$ is needed to say that the contribution of the nonlocal term is still zero since a symmetrization argument is used in Lemma~\ref{first_moment_bdd}.

\section{Extensions -- Hydrodynamic system with alignment}\label{Sec:5}
In this section we shortly discuss how to extend the above results to treat alignment terms at the hydrodynamic level as in \cite{CFGS}. In fact, given a symmetric kernel $\psi\ge0$ with  $\psi\in  L^\infty(\R^3)$, we can consider the following nonlinear damping term in the right hand side of \eqref{main_system} 
\begin{equation}\label{xxx}
\varrho(t,x) \int_{\Omega} \psi (x-y) \Big( \vu(t,y) - \vu(t,x) \Big) \varrho(t,y) \ {\rm d}y 
\end{equation}
instead of the damping term $-\varrho \vu$. In this case, the energy identity \eqref{endis} is changed by replacing the term $\varrho |\vu|^2$ on the right-hand side by 
\begin{equation}\label{xxxx}
\int_\Omega \int_\Omega \psi (x - y) \varrho(t,x) \varrho(t,y) \left| \vu(t,y) - \vu(t,x) \right|^2 \ {\rm d}x {\rm d}y.
\end{equation}
It turns out that for such system there are stationary travelling wave solutions (also called flock solutions) corresponding to constant in space velocity field $\vu=\vu_0\in \R^d$ and $\varrho(t,x)=\varrho_s\left(x-t\vu_0\right)$,
with $\varrho_s$ verifying the balance of forces relation \eqref{balanceforces}. 

The proof of  Theorem \ref{t_existence} can be generalized to this case, the changes are minor and consist in treating the new term \eqref{xxx} as we have done for the interaction forces term $\varrho(\nabla K\ast\varrho)$. This is clear as we may split the term \eqref{xxx} as follows
	\begin{equation}\label{xxxxx}
	\varrho ( \psi \ast (\vu \varrho)) + \varrho \vu (\psi \ast \varrho).
	\end{equation}
Then, due to energy estimates and since $\psi \in L^\infty(\R^3)$, we are able to prove that for each level of approximation  the sequence related to the term \eqref{xxx} converges weakly in $L^q$ for some $q \geq 1$. In consequence, as soon as we know that the density function converges a.e. and strongly in $L^m$ and the momentum converges in $C_{weak}([0,T) ;L^{\frac{2m}{m+1}})$, we are able to characterise properly the  limits for both terms of \eqref{xxxxx}.

 Concerning the long-time asymptotic, solutions are expected to converge generically towards travelling wave flocks in the co-moving frame for the particular case $\Phi =0$ and $\Omega = \R^3$. 
More precisely this happens if $\psi\geq \psi_0>0$,  the solutions to the system 
  	\begin{equation}\label{system_new_variable}
	\begin{split}
	\partial_t 
\varrho  + \Div( \varrho \vu) 
	 = &\, 0 \, \\
	\partial_t (\varrho \vu) + \Div ( \varrho \vu \otimes \vu ) + a \Nabla  \varrho^m 
	= &\, \mu \Delta \vu + \Nabla (\lambda + \mu) \Div \vu- (\Nabla K \ast \varrho) \varrho 
		\\ & -
	\varrho\int_{\R^3} \psi (x- y) \left(\vu (t, y)
	 - \vu(t,x)  \right)\varrho(t,y) \,{\rm{d}}y \,.
	\end{split}
	\end{equation}	
given by the Definition~\ref{defsol}  exist globally and if they satisfy 
\eq{\label{iden}
\frac{{\rm d^2}}{{\rm d}t^2} \int_{\R^3} x \varrho \dx = \frac{{\rm d}}{{\rm d}t} \int_{\R^3} \varrho \vu \dx = 0\,.
}
Note that our existence result, Theorem \ref{t_existence}, cannot be applied in this setting, since there is no confinement, $\Phi = 0$ in $\R^3$.
 Below we  give a sketch of the arguments implying what the long-time asymptotic of such solutions would be.
 On one hand, the dissipative term \eqref{xxxx} provides that $\vu_s(t,x)=\vu_s(t,y)$ for $x,y\in\mbox{supp} (\varrho_s)$. On the other hand, thanks to  identity \eqref{iden} we obtain 
$$
\int_{\R^3} \varrho(t) \vu(t) \dx = \int_{\R^3} \varrho_0 \vu_0 \dx = \int_{\R^3}  \vm_0 \dx
\quad\mbox{ and }\quad 
 \int_{\R^3} x \varrho  (t) \dx = \int_{\R^3} x \varrho_0 \dx + t \int_{\R^3}  \vm_0 \dx
\, .
$$
Then let us define a vector  $\vu_\infty=(u^1_\infty,u^2_\infty,u^3_\infty)$ in $\R^3$ as follows
	$$ M_0 u^i_\infty:=\int_{\R^3} m^i_0 \dx  \quad \mbox{ where } M_0 \mbox{ is the initial mass. }$$ 
Changing variables such that 
$$y=x - \vu_\infty,\quad \tau = t,$$ 
$$\tilde{\varrho}(\tau,y) = \varrho(t, x - t\vu_\infty), \quad\tilde{\vu}(\tau,y)= \vu(t, x-t\vu_\infty),$$ 
we can check that the new velocity and density functions satisfy the very same system of equations \eqref{system_new_variable}. Moreover,  this change of variables leads to 
$$\int_{\R^3} \tilde\varrho \tilde\vu \dy = 0.$$
Note that since $\psi\geq \psi_0>0$, the  integral \eqref{xxxx} implies analogue of \eqref{Lem_LT_1}, and  so, we can proceed as in Section~\ref{s:long_time}.
Therefore, the limit system for the  long time asymptotics appears to be the same as for \eqref{main_system}. Namely, the $\omega$-limit set $\omega(\tilde\varrho)$ associated to the global weak solutions $(\tilde\varrho(\tilde{t}), \, \tilde\vu(\tilde{t}) )$ of \eqref{system_new_variable} consists of stationary solutions with zero momentum and density $\varrho_s$ corresponding to the initial mass $M_0$.   In particular, this solution satisfies the balance of forces relation 
\begin{equation*}\label{balanceforces-xxx}
a \Nabla \varrho^m_s + (\Nabla K \ast \varrho_s) \varrho_s =0\,.
\end{equation*}

%%%%%%%%%%%%%%%%%%%%%%%%%%%%%%%%%%%%%%

\section*{Acknowledgments}
JAC was partially supported by the Royal Society via a Wolfson Research Merit Award and by EPSRC grant number EP/P031587/1. AWK is partially supported by a Newton Fellowship of the Royal Society and by the grant Iuventus Plus no. 0871/IP3/2016/74 of Ministry of Sciences and Higher Education RP. EZ was supported by the UCL Department of Mathematics Grant and grant Iuventus Plus  no. 0888/IP3/2016/74 of Ministry of Sciences and Higher Education RP.

%%%%%%%%%%%%%%%%%%%%%%%%%
%%%%%%%%%%%%%%%%%%%%%%%%%%

\def\ocirc#1{\ifmmode\setbox0=\hbox{$#1$}\dimen0=\ht0 \advance\dimen0
  by1pt\rlap{\hbox to\wd0{\hss\raise\dimen0
  \hbox{\hskip.2em$\scriptscriptstyle\circ$}\hss}}#1\else {\accent"17 #1}\fi}

\end{document}